\RequirePackage{fix-cm}
\documentclass[smallextended]{svjour3}
\smartqed

\newtheorem{thm}{Theorem}
\newtheorem{lem}[thm]{Lemma}
\newtheorem{defn}[thm]{Definition}
\newtheorem{rmk}[thm]{Remark}
\newtheorem{coro}[thm]{Corollary}
\newtheorem{examp}[thm]{Example}
\newtheorem{prop}[thm]{Proposition}

\newcommand{\sk}{\smallskip}

\newcommand{\bt}{\begin{thm}}
\newcommand{\et}{\end{thm}}
\newcommand{\bl}{\begin{lem}}
\newcommand{\el}{\end{lem}}
\newcommand{\bp}{\begin{prop}}
\newcommand{\ep}{\end{prop}}
\newcommand{\br}{\begin{rmk}}
\newcommand{\er}{\end{rmk}}
\newcommand{\bc}{\begin{coro}}
\newcommand{\ec}{\end{coro}}
\newcommand{\bd}{\begin{defn}}
\newcommand{\ed}{\end{defn}}
\newcommand{\bex}{\begin{examp}}
\newcommand{\eex}{\end{examp}}
\newcommand{\beq}{\begin{equation}}
\newcommand{\eeq}{\end{equation}}
\newcommand{\ink}{\hfill $\square $}   

\newcommand{\bce}{\begin{center}}
\newcommand{\ece}{\end{center}}
\newcommand{\pf}{\textit{Proof. }}

\usepackage{graphicx}
\usepackage{amssymb}
\usepackage{graphics, epic, eepic}

\begin{document}

\title{Involutive right-residuated l-groupoids%
\thanks{The work of the first author is supported by the Project I
1923-N25 by Austrian Science Found (FWF), and Czech Grant Agency
(GACR)}}

\titlerunning{Residuated l-groupoids}

\author{Ivan~Chajda \and S\'{a}ndor~Radeleczki }

\institute{Ivan Chajda \at
              Department of Algebra and Geometry, Palacky University \\
              17. listopadu 12, 771 46 Olomouc, Czech Republic \\
              \email{ivan.chajda@upol.cz}           
           \and
           S\'{a}ndor Radeleczki \at
              Mathematical Institute, University of Miskolc \\
              H-3515 Miskolc-Egyetemv\'{a}ros, Hungary \\
              \email{matradi@uni-miskolc.hu}}

\date{}

\maketitle

\begin{abstract}
A common generalization  of orthomodular lattices and residuated
lattices is provided corresponding to bounded lattices with an
involution and sectionally extensive mappings. It turns out that
such a generalization can be based on integral right-residuated
l-groupoids. This general framework is applied to MV-algebras,
orthomodular lattices, Nelson algebras, basic algebras and Heyting
algebras.
\end{abstract}

\section{Introduction}

Residuated lattices were introduced in [17], and they are used in
several branches of mathematics, including areas of ideal lattices
of rings, lattice-ordered groups, formal languages and multi-valued
logic. Right-residuated l-groupoids constitute a natural
generalization of residuated lattices (see e.g. [3]), and their
applications cover even a wider field. We will show, that they
provide a useful framework for propositional calculus in
constructive logic and certain logics related to quantum mechanics,
and some computations in universal algebra.

For instance, let $\mathcal{A}=(A,F)$ be an algebra from a congruence
modular variety, and $[\varphi ,\theta ]$ the commutator of two congruences $%
\varphi ,\theta $. Denote by $0_{A}$ and $1_{A}$ the least and the
greatest element of the congruence lattice $($Con$\mathcal{A},\vee
,\wedge )$, respectively.
In [16], a binary operation $\rightarrow $ on Con$%
\mathcal{A}$ was defined as by the formula:
\begin{center}
$ \alpha \rightarrow \beta :=\bigvee {\{\theta \in \rm{Con}
\mathcal{A}\mid [\alpha ,\theta ]\leq \beta \}}. $
\end{center}

\noindent If the identity $[1_{A},\theta ]=\theta $ holds in
Con$\mathcal{A}$ then, in view of [16], $($Con$\mathcal{A},\vee
,\wedge ,[,],$ $\rightarrow ,0_{A},1_{A})$ is an integral
commutative right-residuated l-groupoid.
\medskip

Although we will not study the consequences of the previous example
in the theory of residuated structures, we can see that integral
commutative right-residuated l-groupoids are not exceptional
structures in algebra, and hence we will investigate the connections
between these structures and lattices having an antitone involution
and so-called sectionally extensive antitone mappings.

 In our paper we study some particular classes of right-residuted
l-groupoids. We aim to show the relevance of these classes of
algebras in several research fields. The paper is structured as
follows. In Section 2 some general notions and facts concerning
right-residuated l-groupoids are presented. In Section 3 we prove
that there is a one-to-one correspondence between involution
lattices with sectionally extensive antitone mappings,\ and
involutive right-residuated l-groupoids satisfying a certain
identity. The case when these residuated l-groupoids form residuated
lattices is characterized. In Section 4 some examples of
right-residuated l-groupoids belonging to the mentioned class are
provided. For instance, we show that residuated lattices
corresponding to Nelson algebras belong to this class. We prove that
sectionally pseudocomplemented lattices admitting an antitone
involution can be characterized as right-residuated l-groupoids
satisfying certain identities. A special attention is paid to those
right-residuated l-groupoids which are defined by lattices with
sectionally antitone involutions. In Section 5 is proved that these
algebras are term equivalent to the so-called basic algebras which
can be viewed as a common generalization of MV-algebras and
orthomodular lattices. The fact that these algebras can be
reconstructed from their implication reduct is shown in Section 6.
Finally, in Section 7, some congruence properties of
right-residuated l-groupoids are investigated.

\section{Preliminaries}

{\bf Definition 1.} By a {\it{right-residuated l-groupoid}} is meant
an algebra

\noindent $\mathcal{G}=(L,\vee ,\wedge ,\odot ,\rightarrow ,0,1)$ of type
(2,2,2,2,0,0) such that

\begin{itemize}
\item[(i)] $(L,\vee,\wedge)$ is a lattice with least element $0$ and
greatest element $1$,

\item[(ii)] $(L,\odot)$ is a groupoid, and $1\odot x=x$, for all $x\in
L $.

\item[(iii)] {$\mathcal{G}$ satisfies the \textit{right-adjointness} property,
that is

$ x\odot y\leq z $ if and only if $x\leq y\rightarrow z$, for all
$x,y,z\in L$ (see e.g. [2]).}

\end{itemize}

\medskip

\noindent In general, right-adjointness does not imply
left-adjointness (see
[4]), except the case when $\mathcal{G}$ is \textit{commutative}, that is, $%
x\odot y=y\odot x$, for all $x,y\in L$.

For our sake, we modify the concept of an integral residuated
structure as follows. The algebra $\mathcal{G}$ will be called
\textit{integral} if $1\odot x=x\odot 1=x $ holds for all $x\in L$.
Clearly, $\mathcal{G}$ is integral whenever it is commutative. Let
$\rceil x:=x\rightarrow 0$. The algebra $\mathcal{G}$
is called \textit{involutive} whenever the mapping $%
x\mapsto \rceil x$, $x \in L $ is an \textit{antitone}
\textit{involution} on $L$, i.e. if $x\leq y$ implies $\rceil y\leq
\rceil x$ and
\begin{center}
\hfill  $ \rceil (\rceil x) = x, $ \hfill ($*$)
\end{center}
for all $ x,y \in L$. The identity ($*$) is called the double
negation law.  Of course, every involutive algebra $\mathcal{G}$
satisfies the double negation law, but not conversely. However, if
$\mathcal{G}$ is a residuated lattice, that is,  $\odot$ is
associative and commutative, then $\mathcal{G}$ is involutive if and
only if it satisfies the double negation law. This is because then
$\mathcal{G}$ satisfies the implication

$$x\leq y \mbox{ implies } y \rightarrow z \leq x \rightarrow z, $$

\noindent for any $ x,y,z \in L$, thus also $ \rceil y = y
\rightarrow 0 \leq x\rightarrow 0 = \rceil x $, for all $ x,y \in
L$, $x\leq y $. Further, we say that $\mathcal{G}$ satisfies
\textit{divisibility} if
$$ (x\rightarrow y)\odot x=x\wedge y, $$

\noindent for every $x,y\in L$. Finally, $\mathcal{G}$ satisfies \textit{%
condition} (C) if
$$ z\leq x\odot y \mbox{ if and only if } y\rightarrow
\rceil x\leq \rceil z, $$ for all $ x,y, z \in L$. The basic
properties of right-residuated l-groupoids are collected in the
following lemma.

\sk

\bl Let $\mathcal{G}=(L,\vee ,\wedge ,\odot ,\rightarrow ,0,1)$ be a
right-residuated l-groupoid. Then

\begin{itemize}

\item[\emph{(i)}] $\rceil 0=1$;

\item[\emph{(ii)}] $a\leq b$ if and only if $a\rightarrow b=1$;

\item[\emph{(iii)}] $a\odot 0=0\odot a=0$, for all $a \in L $;

\item[\emph{(iv)}] $y\leq z$ implies $y\odot x\leq z\odot x$ and $%
x\rightarrow y\leq x\rightarrow z$, for all $ x,y, z \in L$;

\item[\emph{(v)}] $x\odot y\leq y$ and $y\rightarrow z=y\rightarrow (y\wedge
z)$, for all $ x,y, z \in L$;

\item[\emph{(vi)}] if $\mathcal{G}$ satisfies the double negation law
then $\rceil 1=0$.
\end{itemize}

\el

\sk

\pf (i) Since $1\odot 0=0$, we have $1\leq 0\rightarrow 0$, and
hence $1=0\rightarrow 0=\rceil 0$.

\noindent (ii) If $a\leq b$ then $1\odot a=a\leq b$, thus $1\leq
a\rightarrow b$ giving $a\rightarrow b=1$. If $a\rightarrow b=1$, then $%
(a\rightarrow b)\odot a\leq b$ implies $a=1\odot a\leq b$.

\noindent (iii) $a\leq 1=0\rightarrow 0$ yields $a\odot 0=0$, and
$0\leq a\rightarrow 0$ gives $0\odot a=0$.

\noindent (iv) Assume $y\leq z$. Since for all $a,b\in L$, $a\odot
b=a\odot b $ yields
$$ a\leq b\rightarrow (a\odot b),
$$

\noindent we get $y\leq z\leq x\rightarrow (z\odot x)$, whence $y\odot x\leq
z\odot x$.

\noindent Further, $x\rightarrow y\leq x\rightarrow y$ yields
$(x\rightarrow y)\odot
x\leq y\leq z$, whence we deduce $x\rightarrow y\leq x\rightarrow z$, for all $%
x,y,z\in L$.

\noindent (v) Since $x\leq 1=y\rightarrow y$, we obtain $x\odot
y\leq y$, for all $x,y\in L$. Thus $x\odot y\leq z $ if and only if
$ x\odot y\leq y\wedge z$, whence we get $x\leq y\rightarrow z $ if
and only if $ x\leq y\rightarrow \left( y\wedge z\right) $. This
implies $y\rightarrow z=y\rightarrow (y\wedge z)$.

\noindent (vi) The double negation law and (i) imply: $ \rceil 1=
\rceil(\rceil0) = 0$. \ink
 \medskip

 An interrelation between condition (C) and the involutive property is
 stated in the following

 \medskip

\noindent \textbf{Proposition 1} \emph{Let} $\mathcal{G}=(L,\vee
,\wedge
,\odot ,\rightarrow ,0,1)$ \emph{be a right-residuated l-groupoid. Then} $%
\mathcal{G}$ \emph{satisfies the double negation law and condition}
(C)  \emph{if and only if $\mathcal{G}$ is involutive and} $x\odot
y=\rceil (y\rightarrow \rceil x)$ \emph{holds for all} $x,y,z\in L$.

\bigskip

\pf The double negation law yields $\rceil \left( \rceil x\right)
=(x\rightarrow 0)\rightarrow 0=x$. If $x\leq y$ then $x\leq 1\odot
y$, and so by (C) we get $\ \rceil y=y\rightarrow 0=y\rightarrow
\rceil 1\leq \rceil x$. Hence $\mathcal{G}$ is involutive, and (C)
implies $x\odot y\geq z $ if and only if $ y\rightarrow \rceil x\leq
\rceil z $ if and only if $ \rceil (y\rightarrow \rceil x)\geq
\rceil \left( \rceil z\right) =z$. Then $\rceil (y\rightarrow \rceil
x)\geq x\odot y$, and $x\odot y\geq \rceil (y\rightarrow \rceil x)$,
whence $x\odot y=\rceil (y\rightarrow \rceil x)$.

Conversely, suppose that $\mathcal{G}$ is involutive, and $x\odot
y=\rceil (y\rightarrow \rceil x)$ holds. Then clearly, $\mathcal{G}$
satisfies the double negation law, and $\rceil (y\rightarrow \rceil
x)\geq z $ if and only if $ y\rightarrow \rceil x\leq \rceil z$.
This means that $z\leq x\odot y$ if and only if $ y\rightarrow
\rceil x\leq \rceil z$, i.e. (C) holds. \ink

\bigskip

\noindent \textbf{Remark 1} Observe that in a right-residuated l-groupoid the operations $%
\odot $ and $\rightarrow $ determine completely each other, in other words,
 if $\mathcal{G}_{1}=(L,\vee ,\wedge ,\odot ,\rightarrow ,0,1)$ and $\mathcal{%
G}_{2}$ = \ $(L,\vee ,\wedge ,\otimes ,\leadsto ,0,1)$ are
right-residuated l-groupoids having the same underlying lattice
$(L,\vee ,\wedge )$, then the operations $\odot $ and $\otimes $
coincide if and only if $\rightarrow $ and $\leadsto $ coincide. The
proof is the same as that for residuated lattices and hence it is
omitted.

\medskip

Let $\mathcal{G}=(L,\vee ,\wedge ,\odot ,\rightarrow ,0,1)$ be a
right-residuated l-groupoid and define a binary operation $\Rightarrow $ on $%
L $ as follows:%
\begin{center}
$ x\Rightarrow y:=\rceil y\rightarrow \rceil x$, for all $x,y\in L
$.
\end{center}

\noindent Then $\Rightarrow $ will be called the \textit{derived
implication} \textit{of} $\mathcal{G}$.

\medskip

\noindent \textbf{Lemma 2} \emph{Let} $\mathcal{G}=(L,\vee ,\wedge
,\odot ,\rightarrow ,0,1)$ \emph{be an involutive right-residuated
l-groupoid. Then the operation\ }$\Rightarrow $ \emph{for all
$x,y,z\in L$ satisfies the conditions:}

\begin{description}
\item[(I0)] $(x\vee y)\Rightarrow y=x\Rightarrow y$, $x\Rightarrow x=1$, $%
1\Rightarrow x=x$;

\item[(I1)] $(x\Rightarrow y)\wedge y=y$;

\item[(I2)] $x\leq y$ \emph{implies} $y\Rightarrow z\leq x\Rightarrow
z$;
\end{description}

\noindent \emph{Moreover,\ we have} $x\leq y$ \emph{if and only if} $%
x\Rightarrow y=1$.

\bigskip

\pf Since $\mathcal{G}$ is involutive, we have $\rceil 1=0 $, and
hence
\begin{equation}
1\Rightarrow x =\rceil (\rceil x), \mbox{for all } x\in L .
\end{equation}

By definition $x\Rightarrow x=\rceil x\rightarrow \rceil x=1$, and
$(x\vee y)\Rightarrow y=\rceil y\rightarrow \rceil (x\vee y) $, for
all $x,y \in L$. Since $x\mapsto \rceil x$, $ x\in L $ is an
antitone involution on $L$, we have $\rceil (x\vee y)=\rceil x\wedge
\rceil y$, and hence $(x\vee y)\Rightarrow y=\rceil y\rightarrow
(\rceil y\wedge \rceil x)=\rceil y\rightarrow \rceil x=x\Rightarrow
y$, by (v) of Lemma 1. Since $\mathcal{G}$ is involutive, it
satisfies the double negation law, and because ($1$) holds true,
(I0) is clear. By Lemma 1(iv) for any $x,y\in L$ we get $y=\rceil
\left( \rceil y\right) =\rceil y\rightarrow 0\leq \rceil
y\rightarrow \rceil x=x\Rightarrow y$, which proves (I1).

\noindent (I2). Since $\mathcal{G}$ is involutive, we have $x\leq y
$ if and only if $ \rceil y\leq \rceil x$. By Lemma 1(iv) $\rceil
y\leq \rceil x$ implies $\rceil z\rightarrow \rceil y\leq \rceil
z\rightarrow \rceil x$. Hence $x\leq y$ implies $y\Rightarrow z\leq
x\Rightarrow z$.

\noindent Finally, $x\leq y $ if and only if $ \rceil y\leq \rceil x
$, and Lemma 1(ii) yields $ \rceil y\leq \rceil x $ if and only if
$\rceil y\rightarrow \rceil x=1 $. However $\rceil y\rightarrow
\rceil x=1 $ means that $ x\Rightarrow y=1$. \ink

\medskip

\section{Lattices with sectionally antitone mappings}

\bigskip

 An algebraic axiomatization of \L ukasie\-wicz many-valued logic
 can be
provided by means of MV-algebras, and analogously, orthomodular
lattices constitute an important algebraic framework for logical
computations related to quantum mechanics. As will be shown in
Section 4, both of these classes of algebras can be recognized as
bounded lattices with sectionally antitone involutions. However, not
in all the algebraic structures used for the formalization of
non-classical logics the corresponding sectional mappings (derived
by the logical connective implication) must be involutions. For
example, in the case of Heyting algebras or BCK-algebras these
mappings are antitone, but not necessarily they are involutions.
Hence we introduce formally the concept of a lattice with
sectionally antitone mappings which will be used here.

\medskip

\noindent Let $(L,\vee ,\wedge ,0,1)$ be a bounded lattice. For an
$a\in L$ the interval $[a,1]=\{x\in L\\\mid a\leq x\leq 1\}$ is called
a \textit{section}. The algebra $\mathcal{L}=(L,\vee ,\wedge
,\{^{a}\mid a\in L\},0,1)$ is called \textit{a lattice with
sectionally antitone extensive mappings }if for each $a\in L$ there
exists a mapping $x\mapsto x^{a}$ of $[a,1]$ into itself, such that
\begin{center}
 \ \ \ $ x\leq y $ implies  $ x^{a}\geq y^{a}$, for all $x,y\in \lbrack
a,1] $, and \hfill (i.e. $x\mapsto x^{a}$ is antitone)
\end{center}

\begin{center}
 \ \ \  $ x^{aa}\geq x$, for all $x\in \lbrack a,1]$. \hfill (i.e.
$x\mapsto x^{a}$ is extensive)
\end{center}

\noindent In this case $1^{a}=a$ implies $a^{a}=1$. Indeed, $1^{aa}=1$
yields $a^{a}=(1^{a})^{a}=1$.

In particular, if each mapping $x\mapsto x^{a}$, $x \in [a,1] $ is an involution, i.e. $%
x^{aa}=x$, for all $x \in [a,1] $, then $\mathcal{L}$ is called a
\textit{lattice with sectionally antitone involutions }(see e.g.
[8]).

\medskip

Let us note that in our example $($Con$\mathcal{A},\vee ,\wedge
,[,],$ $\rightarrow ,0_{A},1_{A})$ from the introduction, for
any $ \alpha, \theta \in $ Con$%
\mathcal{A}$, with  $ \alpha \leq \theta $ we can define

\begin{center}
 $\theta^{\alpha}:= \theta \rightarrow \alpha=\bigvee {\{\varphi \in \rm{Con} \mathcal{A}\mid
[\theta ,\varphi ]\leq \alpha \}}. $
\end{center}

\noindent Since $ [\theta ,\varphi ]\leq \theta \wedge \varphi $
holds in any congruence modular variety, we get $[\theta ,\alpha
]\leq \alpha $, and hence $\theta^{\alpha} \geq \alpha $. Since for
any $ \theta_{1}, \theta_{2}, \varphi \in \rm{Con} \mathcal{A}$
$\theta_{1} \leq \theta_{2} $ implies $ [\theta_{1} ,\varphi ]\leq
[\theta_{2} ,\varphi ]$, we get $\theta_{1}^{\alpha} \geq
\theta_{2}^{\alpha}$ whenever $ \alpha \leq \theta_{1} \leq
\theta_{2} $. Finally, $ [\theta^{\alpha},\theta]=[\theta
\rightarrow \alpha, \theta]\leq \alpha $ implies $ \theta^{\alpha
\alpha}\geq \theta $. Thus for any $ \alpha \in \rm{Con}
\mathcal{A}$ the mapping $ \theta \mapsto \theta^{\alpha} $, $\theta
\in [\alpha,1_{A}]$ is a sectionally antitone extensive mapping.
\medskip

\noindent \textbf{Proposition 2} \emph{Let }$(L,\vee ,\wedge
,0,1)$\emph{\
be a bounded lattice and }$\Rightarrow $ \emph{a} \emph{binary operation on }%
$L$,  \emph{and define} $ x^{a}:=x\Rightarrow a$, \emph{for any }$a,
x\in L$, \emph{with }$x \geq a$. \emph{Then the following are
equivalent:}

\begin{itemize}
\item[(i)] \emph{The binary operation }$\Rightarrow $ \emph{satisfies }(I0),
(I1), (I2)\emph{\ and}

\begin{center}
$ [(x\Rightarrow y)\Rightarrow y] \wedge (x\vee y)=(x\vee y) $, for
all  $x,y\in L $. \hfill (I3)
\end{center}

\item[(ii)] \emph{For each} $ a \in L$ \emph{the mapping }$x\mapsto x^{a}$, $x \in [a,1] $,
\emph{\ is an antitone extensive mapping on} $[a,1]$ \emph{such
that} $1^{a}=a$ \emph{and} $x\Rightarrow y=(x\vee y)^{y}$\emph{,}
\emph{for all }$x,y\in L$\emph{.}
\end{itemize}

\bigskip

\pf (i)$\Rightarrow $(ii). Take $a,x\in L$ arbitrary
with $x\geq a$. Then in view of (I1) we get $%
a\leq x\Rightarrow a=x^{a}$, and this means that the assignment
$x\mapsto x^{a}$, $x\in \lbrack a,1]$ is a mapping of $[a,1]$ into
itself. Let $a\leq x\leq y$. Then (I2) yields $y^{a}=y\Rightarrow
a\leq x\Rightarrow a=x^{a}$, hence $x\mapsto
x^{a}$, $x\in \lbrack a,1]$ is antitone. By using (I3), for every $x\in \lbrack a,1]$ we obtain $%
x^{aa}=(x\Rightarrow a)\Rightarrow a\geq x\vee a=x$, i.e. the mapping $%
x\mapsto x^{a}$, $x\in \lbrack a,1]$ is extensive. Finally, (I0) implies $1^{a}=1\Rightarrow a=a$%
, and $x\Rightarrow y=(x\vee y)\Rightarrow y=(x\vee y)^{y}$, for all
$x,y\in L $.

(ii)$\Rightarrow $(i). Let $\mathcal{L}=(L,\vee ,\wedge ,\{^{a}\mid a\in
L\},0,1)$ be a lattice with sectionally extensive antitone mappings $%
x\mapsto x^{a}$, $x\in \lbrack a,1]$ such that $1^{a}=a$, for all
$a\in L$, and suppose that, for
all  $x,y\in L $ the operation $\Rightarrow $ satisfies%
\[
x\Rightarrow y=(x\vee y)^{y}.
\]

\noindent Then $(x\vee y)\Rightarrow y=x\Rightarrow y$. Since
$1^{a}=a$ implies $a^{a}=1$, we get $x\Rightarrow x=x^{x}=1$ and
$x\Rightarrow 1=1\Rightarrow 1=1$, and also $1\Rightarrow
x=1^{x}=x$, for all $ x \in L $. Thus (I0) is satisfied. As by
definition $x\Rightarrow y=(x\vee y)^{y}\geq y$, we get
$(x\Rightarrow y)\wedge y=y$, for all $x,y\in L$, i.e. (I1) holds.
Now assume $x\leq y$. Then $x\vee z\leq y\vee z$, for all $z\in L$%
, and hence $y\Rightarrow z=(y\vee z)^{z}\leq (x\vee
z)^{z}=x\Rightarrow z$, for all $x,y,z\in L$ because the map
$x\mapsto x^{z}$, $x\in \lbrack z,1]$ is antitone. Thus (I2) holds
for $\Rightarrow $. To prove (I3), let us observe that
$(x\Rightarrow y)\Rightarrow y=(\left( x\vee y\right) \Rightarrow
y)\Rightarrow y=(x\vee y)^{yy}$, for all  $x,y\in L $. Since by
extensive property $(x\vee y)^{yy}\geq x\vee y$, we obtain $\left[
(x\Rightarrow y)\Rightarrow y\right] \wedge (x\vee y)=(x\vee
y)^{yy}\wedge (x\vee y)=x\vee y$, for all $x,y\in L$. \ink

\medskip

The mutual interrelation between involutive right-residuated
l-groupoids satisfying condition (I3) and bounded lattices with an
antitone involution and sectionally extensive antitone mappings is
established in the next theorem. This gives us an alternative
approach to involutive right-residuated l-groupoids which is more
suitable to algebras used for axiomatization of several
non-classical logics.

\medskip

\noindent\textbf{Theorem 1}
\begin{itemize}
\item[(a)] \emph{Let} $\mathcal{L}=(L,\vee,\wedge ,\{^{a}\mid a\in
L\},\thicksim,0,1)$ \emph{be a bounded lattice with an antitone
involution }$\thicksim$ \emph{and sectionally antitone extensive
mappings }$x\mapsto x^{a},$\emph{\ }$x\in\lbrack a,1]$\emph{\ such that }$%
1^{a}=a,$\emph{\ for all }$a\in L$\emph{. If we define }

\begin{equation}
x\rightarrow y:=(\thicksim x\vee\thicksim y)^{\thicksim x}
\end{equation}

\begin{equation}
x\odot y:=\thicksim(y\rightarrow\thicksim
x)=\thicksim\lbrack(x\vee\thicksim y)^{\thicksim y}],
\end{equation}

\noindent \emph{for all} $x,y \in L$\emph{, then}
$\mathcal{G}(\mathcal{L})=(L,\vee ,\wedge ,\odot ,\rightarrow ,0,1)$
\emph{is an involutive right-residuated l-groupoid such
that }$\rceil x=\;\thicksim x$\emph{ holds, and its derived implication }$%
x\Rightarrow y:=\rceil y\rightarrow \rceil x$\emph{\ satisfies
condition }(I3)\emph{.}

\item[(b)] \emph{Let} $\mathcal{G}=(L,\vee ,\wedge ,\odot ,\rightarrow ,0,1)$
\emph{be an involutive right-residuated l-groupoid having the
property that its derived implication }$\Rightarrow $\emph{\
satisfies condition }(I3). \emph{Let }$\thicksim z:=z\rightarrow
0$\emph{, for all }$z\in L$\emph{, and
define }%
\begin{equation}
x^{a}:=x\Rightarrow a=\rceil a\rightarrow \rceil x,
\end{equation}

\noindent \emph{for all }$a,x\in L$\emph{\ with }$x\geq a$\emph{. Then} $%
\mathcal{L}(\mathcal{G})=(L,\vee ,\wedge ,\{^{a}\mid a\in L\},\thicksim
,0,1) $ \emph{is a bounded lattice with an antitone involution }$\thicksim $
\emph{and sectionally antitone extensive mappings }$x\mapsto x^{a},$\emph{\ }%
$x\in \lbrack a,1]$\emph{\ such that }$1^{a}=a$\emph{.}

\item[(c)] \emph{The correspondence between bounded lattices with an involution }$\thicksim $ \emph{ and sectionally
antitone extensive mappings satisfying} $1^{a}=a$, \emph{and
involutive right-residuated l-groupoids
satisfying condition }(I3)\emph{\ is one-to-one, i.e.} $\mathcal{G}(\mathcal{L}%
(\mathcal{G}))=\mathcal{G}$ \emph{and} $\mathcal{L}(\mathcal{G}(\mathcal{L}%
))=\mathcal{L}$.
\end{itemize}
\medskip

 Before the proof, let us note that the mappings $x\mapsto $ $%
\sim $ $x$,  $x \in L$ and $x\mapsto x^{0}$, $x \in L$ need not
coincide. The second map
need not be an involution contrary to the case $%
x\mapsto $ $\sim $ $x$, $x \in L$.

\medskip

\pf (a) By definition we have
\begin{center}
$1\odot x=\thicksim\lbrack (1\vee\thicksim x)^{\thicksim
x}]=\thicksim(1^{\thicksim x})=\thicksim (\thicksim x)=x$, for all $
x \in L$. \hfill($\ast$)
\end{center}

\noindent Let $x\odot y\leq z$ for some $x,y,z\in L$. Then
$\thicksim \lbrack(x\vee\thicksim y)^{\thicksim y}]\leq z$ implies
that $\thicksim z\leq (x\vee\thicksim y)^{\thicksim y}$. Since
$\thicksim y\leq(x\vee\thicksim y)^{\thicksim y}$, together we
obtain
\[
\thicksim z\vee\thicksim y\leq(x\vee\thicksim y)^{\thicksim y}.
\]

\noindent This implies $x\leq x\vee \thicksim y\leq (x\vee \thicksim
y)^{\thicksim y\thicksim y}\leq \left( \thicksim z\vee \thicksim
y\right) ^{\thicksim y}=y\rightarrow z$, according to the definition
and to the antitony of the mapping $x\mapsto x^{\thicksim y}$, $x
\in [\thicksim y, 1]$.

Conversely, $x\leq y\rightarrow z$ implies $x\vee \thicksim y\leq
\left( \thicksim z\vee \thicksim y\right) ^{\thicksim y}$, whence we
get $\left( \thicksim z\vee \thicksim y\right) ^{\thicksim
y\thicksim y}\leq (x\vee \thicksim y)^{\thicksim y}$, thus
$\thicksim \lbrack (x\vee \thicksim y)^{\thicksim y}]\leq $
$\thicksim \lbrack \left( \thicksim z\vee \thicksim y\right)
^{\thicksim y\thicksim y}]$. Because the map $x\mapsto x^{\thicksim
y}$, $ x \in [\thicksim y,1] $ is extensive $\left( \thicksim z\vee
\thicksim y\right)
^{\thicksim y\thicksim y}\geq $ $\thicksim z$, whence we deduce $%
\thicksim \lbrack \left( \thicksim z\vee \thicksim y\right) ^{\thicksim
y\thicksim y}]\leq \thicksim (\thicksim z)=z$. Thus we obtain:%
\[
x\odot y=\thicksim \lbrack (x\vee \thicksim y)^{\thicksim y}]\leq z.
\]

\noindent Since $\mathcal{G}(\mathcal{L})$ satisfies the right-adjointness
property and ($\ast$), it is a right residuated l-groupoid. Observe also, that $%
\rceil x:=x\rightarrow 0=(\thicksim x\vee \thicksim 0)^{\thicksim
x}=1^{\thicksim x}=\thicksim x$. Thus the map $x\mapsto \rceil x$, $
x \in L $ is an
antitone involution on $L$, and we can write:%

\[
x\rightarrow y=(\rceil x\vee \rceil y)^{\rceil x}, \ \ x\odot
y=\rceil (y\rightarrow \rceil x)=\rceil (x\vee \rceil y)^{\rceil y},
\] and %
\[
x\Rightarrow y=\rceil y\rightarrow \rceil x=(x\vee y)^{y}.
\]

\noindent Hence for any $a\in L$ and $x\in \lbrack a,1]$ we get
$x^{a}=(x\vee a)^{a}=x\Rightarrow a$. Then $\Rightarrow $ satisfies
(I3), according to Proposition 2.

\noindent (b) Since $\mathcal{G}=(L,\vee ,\wedge ,\odot ,\rightarrow
,0,1)$ is involutive, the map $\thicksim x:=x\rightarrow 0=\rceil
x$, $x\in L$ is an antitone involution, and by using Lemma 2 we get
that $x\Rightarrow y=\rceil y\rightarrow \rceil x$ satisfies
(I0),(I1) and (I2). Since (I3) is also satisfied by $\Rightarrow $,
by defining $x^{a}:=x\Rightarrow a$, for all $a\in L$ and $x\in
\lbrack a,1]$, and using Proposition 2, we obtain that
$\mathcal{L}(\mathcal{G})=(L,\vee ,\wedge ,\{^{a}\mid a\in L\},0,1)$
is
a lattice with sectionally antitone extensive mappings $x\mapsto x^{a},$%
\emph{\ }$x\in \lbrack a,1]$ satisfying $1^{a}=a$.

\noindent (c) First, we prove that $\mathcal{G}(\mathcal{L}(\mathcal{G}))=%
\mathcal{G}$.

Indeed, in $\mathcal{L}(\mathcal{G)}$ we have $\thicksim
x=x\rightarrow 0=\rceil x$, for all $ x \in L$ where $\rceil
x:=x\rightarrow 0$ is defined in $\mathcal{G}$. Then by (a), $\rceil
x$ has the same meaning as in
$\mathcal{G}(\mathcal{L}(\mathcal{G}))$. In view of $(2)$, for all
$x,y\in L$ the operation $\rightarrow $ in
$\mathcal{G}(\mathcal{L}(\mathcal{G}))$ is defined as

$x\rightarrow y:=(\thicksim x\vee \thicksim y)^{\thicksim x}=(\rceil x\vee
\rceil y)^{\rceil x}=(\rceil x\vee \rceil y)\Rightarrow \rceil x$, where $%
\Rightarrow $ is the derived implication of $\mathcal{G}$. Since
(I0) holds in $\mathcal{G}$, we get $(\rceil x\vee \rceil
y)\Rightarrow \rceil x=\rceil y\Rightarrow \rceil x$. Thus we obtain
$x\rightarrow y=\rceil y\Rightarrow \rceil x$. Since in view of (b),
$\rceil y\Rightarrow \rceil x$ also equals to $x\rightarrow y$ in
$\mathcal{G}$, the operation $\rightarrow $ in the right-residuated
l-groupoid $\mathcal{G}(\mathcal{L}(\mathcal{G}))$ coincides with
the operation $\rightarrow $ in $\mathcal{G}$. Therefore, in view of
Remark 1, $\odot $ represents the same operation in $\mathcal{G}$ and $%
\mathcal{G}(\mathcal{L}(\mathcal{G}))$. Because these algebras are defined
on the same bounded lattice $(L,\vee ,\wedge ,0,1)$, they coincide, i.e. $%
\mathcal{G}(\mathcal{L}(\mathcal{G}))=\mathcal{G}$.

To prove $\mathcal{L}(\mathcal{G}(\mathcal{L}))=\mathcal{L}$, first observe
that for any $x\in L$, $\thicksim x$ in $\mathcal{L}(\mathcal{G}(\mathcal{L}%
))$ is defined as $x\rightarrow 0=\rceil x$ in $\mathcal{G}(\mathcal{L})$,
and this is the same as $\thicksim x$ in $\mathcal{L}$, according to (a).
Hence the algebras $\mathcal{L}$ and $\mathcal{L}(\mathcal{G}(\mathcal{L}))$
are defined on the same bounded lattice $(L,\vee ,\wedge ,\thicksim ,0,1)$
with an antitone involution. Therefore, it is enough to prove that the
mappings $x\mapsto x^{a},$\emph{\ }$x\in \lbrack a,1]$ are the same in $%
\mathcal{L}(\mathcal{G}(\mathcal{L}))$ and $\mathcal{L}$. Observe that $%
x^{a} $ in $\mathcal{L}(\mathcal{G}(\mathcal{L}))$ by definition is the same
as $\rceil a\rightarrow \rceil x$ in the right-residuated l-groupoid $%
\mathcal{G}(\mathcal{L})$. By the definition of $\mathcal{G}(\mathcal{L})$
in (a) we get
\[
\rceil a\rightarrow \rceil x=\;\thicksim a\rightarrow \/\thicksim
x=(a\vee x)^{a}=x^{a},
\]

\noindent where $x^{a}$ is defined in $\mathcal{L}$ for all $a,x \in $ with $x\geq a$. Hence $%
x^{a}$ in $\mathcal{L}(\mathcal{G}(\mathcal{L}))$ is the same as $x^{a}$ in $%
\mathcal{L}$, and this completes the proof. \ink

\bigskip

\noindent \textbf{Corollary 1 }\emph{Let} $\mathcal{G}=(L,\vee
,\wedge ,\odot ,\rightarrow ,0,1)$ \emph{be an involutive
right-residuated l-groupoid. Then the following assertions are
equivalent:}

\begin{itemize}
\item[(i)] \emph{The derived implication }$\Rightarrow $\emph{\ satisfies
identity }(I3).

\item[(ii)] $x\odot y=\rceil (y\rightarrow \rceil x)$ \emph{holds for
all} $x,y \in L$.

\item[(iii)] $\mathcal{G}$ \emph{satisfies} \emph{condition} (C).
\end{itemize}

\medskip

\pf Since $\mathcal{G}$ satisfies the double negation law, in view
of Proposition 1, (ii) and (iii) are equivalent.

\noindent (i)$\Rightarrow $(ii). If (i) holds then $\Rightarrow $ satisfies
all the conditions (I0),...,(I3), according to Lemma 2. Now (ii) follows by
applying Proposition 2 and Theorem 1.

\noindent (ii)$\Rightarrow $(i). Since $\mathcal{G}$ is involutive,
in view of Lemma 2, $\Rightarrow $ satisfies (I1). This implies
$y\leq (x\Rightarrow y)\Rightarrow y$, for any $x ,y\in L$. Observe
that in order to prove (I3) it is enough to show that $x\leq
(x\Rightarrow y)\Rightarrow y$. We have:

$(x\Rightarrow y)\Rightarrow y=\rceil y\rightarrow \rceil
(x\Rightarrow y)=\rceil y\rightarrow \rceil (\rceil y\rightarrow
\rceil x)=\rceil y\rightarrow (x\odot \rceil y)$. Now, $x\odot
\rceil y\leq x\odot \rceil y$ gives $x\leq \rceil y\rightarrow
(x\odot \rceil y)=(x\Rightarrow y)\Rightarrow y$, completing the
proof. \ink

\medskip

 Observe that residuated lattices can be characterized as integral
residuated l-groupoids where the operation $\odot $ is associative
and commutative. Hence it is important in our case to know under
what conditions the above properties hold.

\bigskip

\noindent \textbf{Theorem 2} Let $\mathcal{G}=(L,\vee ,\wedge ,\odot
,\rightarrow ,0,1)$ \emph{be an involutive right-residuated
l-groupoid
satisfying }$x\odot y=\rceil (y\rightarrow \rceil x)$ \emph{for all} $x,y \in L$ \emph{and }$%
\Rightarrow $ \emph{its derived implication.} \emph{Then the
following hold true:}

\begin{itemize}

\item[(i)] $\mathcal{G}$ \emph{is integral if and only if} $x\Rightarrow
0=x\rightarrow 0$, \emph{for all }$x\in L$\emph{.}

\item[(ii)] $\mathcal{G}$ \emph{is commutative if and only if }$%
\Rightarrow $ \emph{and }$\rightarrow $\emph{\ coincide.}

\item[(iii)] $\odot $\emph{\ is associative if and only if }
\begin{center}
$ (x\odot y)\Rightarrow z=x\Rightarrow (y\Rightarrow z) $, \emph{for
all} $x,y \in L$. \hfill (D)
\end{center}

\end{itemize}

\medskip

\pf (i) If $1\odot x=x\odot 1=x$ holds for all $x\in L$, then $%
x\leq 1\rightarrow x$, and $1\rightarrow x=(1\rightarrow x)\odot
1\leq x$, hence $x=1\rightarrow x$. Then $x\rightarrow 0=\rceil
x=1\rightarrow \rceil x=\rceil (\rceil x)\Rightarrow \rceil
1=x\Rightarrow 0$, because $\mathcal{G} $ satisfies the double
negation law.

Conversely, suppose that $x\Rightarrow 0=x\rightarrow 0$, for all\emph{\ }$%
x\in L$. Then $x\odot 1=\rceil (1\rightarrow \rceil x)=\rceil (\rceil
(\rceil x)\Rightarrow \rceil 1)=\rceil (x\Rightarrow 0)=\rceil (x\rightarrow
0)=\rceil (\rceil x)=x$.

\noindent (ii) By our assumption, $x\odot \rceil y=\rceil (\rceil
y\rightarrow \rceil x)=\rceil (x\Rightarrow y)$. Hence,
$x\Rightarrow y=\rceil (x\odot \rceil y)$, for all $x,y \in L$. If
$\odot $ is commutative, then $x\Rightarrow y=\rceil (x\odot \rceil
y)=\rceil (\rceil y\odot x)=\rceil (\rceil (x\rightarrow \rceil
(\rceil y))=x\rightarrow y$, for all $x,y \in L$.

Conversely, $x\Rightarrow y=x\rightarrow y$ implies $x\Rightarrow
\rceil y=x\rightarrow \rceil y$. This means that $\rceil (\rceil
y)\rightarrow \rceil x=x\rightarrow \rceil y$, i.e. $y\rightarrow
\rceil x=x\rightarrow \rceil y$. Then for all $x,y \in L$ we have
$x\odot y=\rceil (y\rightarrow \rceil x)=\rceil (x\rightarrow \rceil
y)=y\odot x$, hence $\mathcal{G}$ is commutative.

\noindent (iii) We have $(x\odot y)\odot z=\rceil (z\rightarrow
\rceil (x\odot y))=\rceil (z\rightarrow (y\rightarrow \rceil x))$,
for all $x,y,z\in L$. Observe that $(x\odot y)\Rightarrow \rceil
z=\rceil \left( \rceil
z\right) \rightarrow \rceil (x\odot y)=z\rightarrow (y\rightarrow \rceil x)$%
. Hence $(x\odot y)\odot z=\rceil ((x\odot y)\Rightarrow \rceil z)$.
Similarly, we get $x\odot (y\odot z)=\rceil ((y\odot z)\rightarrow \rceil
x)=\rceil (\rceil (z\rightarrow \rceil y)\rightarrow \rceil x)=\rceil
(x\Rightarrow (z\rightarrow \rceil y))=\rceil (x\Rightarrow (y\Rightarrow
\rceil z))$.

First, suppose that $\odot $ is associative. Then $(x\odot y)\odot z=x\odot
(y\odot z)$ implies%
\[
(x\odot y)\Rightarrow \rceil z=x\Rightarrow (y\Rightarrow \rceil z),
\]

\noindent and $(x\odot y)\Rightarrow z=(x\odot y)\Rightarrow \rceil
(\rceil z)=x\Rightarrow (y\Rightarrow \rceil (\rceil
z))=x\Rightarrow (y\Rightarrow z)$, for all $x,y,z \in L$, which is
(D).

 Conversely, suppose that (D) holds. Then $\rceil \left( (x\odot
y)\Rightarrow \rceil z\right) =\rceil \left( x\Rightarrow
(y\Rightarrow \rceil z)\right) $ is also satisfied, for all $x,y,z
\in L$. In view of the above formulas, this means that $(x\odot
y)\odot z=x\odot (y\odot z)$, for all $x,y,z \in L$. Thus $\odot $
is associative. \ink

\medskip

\noindent \textbf{Corollary 2 }Let $\mathcal{G}=(L,\vee ,\wedge
,\odot ,\rightarrow ,0,1)$ \emph{be an involutive right-residuated
l-groupoid such
that }$\Rightarrow $\ \emph{satisfies condition} (I3)\emph{. Then }$%
\mathcal{G}$ \emph{is an integral commutative residuated lattice if and only
if }$\odot $\emph{\ is associative.}

\medskip

\pf  Since the only if part is clear, and $\mathcal{G}$ is integral
whenever it is commutative, we have to show only that $\odot $ is
commutative, whenever it is associative.

Suppose that $\odot $ is associative. Since we have $x\odot y=\rceil
(y\rightarrow \rceil x)$ by Corollary 1, Theorem 2 yields $(x\odot
y)\Rightarrow z=x\Rightarrow (y\Rightarrow z)$. Then Lemma 2
implies:

\smallskip
\indent \indent $x\odot y\leq z\Leftrightarrow (x\odot y)\Rightarrow
z=1\Leftrightarrow x\Rightarrow (y\Rightarrow z)=1\Leftrightarrow
x\leq y\Rightarrow z$.
\smallskip

\noindent Thus we get $x\leq y\rightarrow z$ if and only if $x\odot
y\leq z$ if and only if $x\leq y\Rightarrow z$, and this implies
$y\rightarrow z\leq y\Rightarrow z$ and $y\Rightarrow z\leq
y\rightarrow z$. Hence $y\rightarrow z=y\Rightarrow z$, for all
$y,z\in L$, and now by using Theorem 2(ii) we obtain that $\odot $
is commutative. \ink

\medskip

It is known that any integral commutative residuated lattice $\mathcal{L}$
satisfying the double negation is involutive (see e.g. [20]). Moreover, $%
x\odot y=\rceil (y\rightarrow \rceil x)$ holds in $\mathcal{L}$, according
to [2; Theorem 2.40]. Hence, by Theorem 3(ii) $\Rightarrow $\emph{\ }and $%
\rightarrow $\ coincide in $\mathcal{L}$, and in view of Corollary 1 and
Theorem 1(b) we obtain:

\medskip

\noindent \textbf{Corollary 3 }\emph{Let
}$\mathcal{L}\emph{=}(L;\vee ,\wedge ,\odot ,\rightarrow ,0,1)$
\emph{be a (commutative, integral)
residuated lattice satisfying the double negation law. Then }$\Rightarrow $%
\emph{\ and }$\rightarrow $\emph{\ coincide, and for each }$a\in L$\emph{, }$%
x^{a}:=x\rightarrow a$, $x\in \lbrack a,1]$\emph{\ is an} \emph{antitone
extensive mapping}.

\bigskip

\section{Examples and applications}

\subsection{Sectionally pseudocomplemented lattices with an added involution}

\medskip

In this section we show how useful can be lattices with an antitone
involution and sectionally extensive mappings. This will be shown by
examples of algebras used frequently in mathematics as well as in
applications.

A bounded lattice $L$ is called
\textit{pseudocomplemented} if for any $x\in L$ there exists an
element $x^{\ast }\in L$ such that
\begin{center}
$y\wedge x=0 $ if and only if $y\leq x^{\ast }$.
\end{center}

\noindent It is evident that $x^{\ast \ast }\geq x$, and $x\leq y$ implies $%
y^{\ast }\leq x^{\ast }$, for any $x,y \in L $. If for any $a\in L$
the section $[a,1]$ is a pseudocomplemented lattice, then $L$ is
called \textit{sectionally pseudocomplemented}.

It is worth mentioning that sectionally pseudocomplemented  lattices
capture the relativity of the pseudocomplement slightly better than
the so-called relatively pseudocomplemented lattices. Namely in a
relatively pseudocomplemented lattice $ L $ the relative
pseudocomplement $x \rightarrow y $ of an element $ x \in L $ with
respect to $ y \in L $ need not to belong to the interval $ [y,1] $,
however it is known that any relatively pseudocomplemented bounded
lattice is also sectionally pseudocomplemented (see [6]). Moreover,
as it is shown in [6], sectionally pseudocomplemented lattices
enable us to extend the concept of relative pseudocomplementation
also for nondistributive lattices. For instance, in [11] is proved
that any algebraic $\wedge$-semidistributive lattice is sectionally
pseudocomplemented; in particular, finite sublattices of free
lattices are sectionally pseudocomplemented lattices which are not
distributive, in general.

\medskip

Let $ L $ be a bounded sectionally pseudocomplemented lattice. For any $ a \in L $ denote by  $x^{a}$ the pseudocomplement of an element $%
x\in \lbrack a,1]$ in the sublattice $([a,1],\leq )$, and define $%
x\Rightarrow y:=(x\vee y)^{y}$, for all  $x,y \in L $. Observe that
$x\mapsto x^{a}$, $x\in \lbrack a,1]$ is an antitone extensive
mapping of $[a,1]$ into itself for each $a\in L$.

Indeed, $x^{a}\in \lbrack a,1]$ by definition, and for any $a\leq
x\leq
y$ we have $y^{a}\leq x^{a}$, and $x^{aa}\geq x$. Then by Proposition 2, $%
\Rightarrow $ satisfies the conditions (I0),..,(I3).

\noindent Now let $\sim $ be an antitone involution on $L$. If we define

$x\rightarrow y:=(\thicksim x\vee \thicksim y)^{\thicksim x}$ and $x\odot
y:=\thicksim \lbrack (x\vee \thicksim y)^{\thicksim y}]=\thicksim
(x\Rightarrow \thicksim y)$,

\noindent for all $x,y \in L $, then by Theorem 1(a) we obtain an
involutive right-residuated
l-groupoid $\mathcal{G}=(L,\vee ,\wedge ,\odot ,\rightarrow ,0,1)$ such that $%
\rceil x=x\rightarrow 0=\thicksim x$, for all $ x \in L$, and its
derived implication
coincides with $\Rightarrow $. \\

A well known example for a sectionally pseudocomplemented lattice
admitting an antitone involution is the five element nondistributive
lattice $N_{5}$. In view of [6] and [11] sectionally
pseudocomplemented bounded lattices are characterized by the
following identities:
\begin{enumerate}

\item[(P1)] $x\Rightarrow x=1$, $1\Rightarrow x=x$, for all $x \in
L $;

\item[(P2)] $(x\vee y)\Rightarrow y=x\Rightarrow y$, $y\wedge
(x\Rightarrow y)=y$, for all $x,y \in L $;

\item[(P3)] $\left[ (x\Rightarrow y)\Rightarrow y\right] \wedge (x\vee
y)=(x\vee y)$, for all $x,y \in L $;

\item[(P4)] $([(x\vee z)\wedge (y\vee z)]\Rightarrow z)\wedge ([(x\vee
z)\wedge (y\Rightarrow z)]\Rightarrow z)=x\wedge z$, for all $x,y,z
\in L $.
\end{enumerate}

\noindent Let us observe that the conjunction of (P1), (P2) and (P3)
is equivalent to the conjunction of (I0), (I1), (I2) and (I3). By
the above characterization $\Rightarrow $ in $\mathcal{G}$ also
satisfies (P4). Moreover, using this characterization and Theorem 1,
we deduce:

\medskip

\noindent \textbf{Proposition 3} \emph{Let} $\mathcal{G}=(L,\vee
,\wedge ,\odot ,\rightarrow ,0,1)$ \emph{be an involutive
right-residuated l-groupoid. Then its derived implication
}$\Rightarrow $ \emph{satisfies condition} (P3) \emph{and} (P4)
\emph{if and only if} $(L,\vee ,\wedge )$ \emph{is a sectionally
pseudocomplemented lattice with an antitone involution such that for
any} $x,y\in L$ \emph{with} $x\geq y$, $x\Rightarrow y$ \emph{is
equal to the pseudocomplement of }$x$\emph{\ in }$[y,1]$\emph{.}

\bigskip

\noindent We note that $\mathcal{G}$ is neither integral nor associative, in
general. Clearly, if $\odot $ is associative, then $\mathcal{G}$ is integral
by Corollary 2. If $\mathcal{G}$ is integral, then we have $x^{\ast
}=x\Rightarrow 0=x\rightarrow 0=$\/$\thicksim x$, according to Theorem 2. It
is known that the map $x\mapsto x^{\ast }$, $x \in L $ is an involution on $L$ if and only if $%
(L,\vee ,\wedge )$ is a Boolean lattice. Hence $\mathcal{G}$ is
integral if and only if $(L,\vee ,\wedge )$ is a Boolean lattice.

\bigskip

\subsection{Residuated lattices corresponding to Nelson algebras}

\medskip

 Let $(L,\vee ,\wedge ,0,1)$ be a bounded distributive lattice with an
antitone involution $\sim $. If for all $x,y\in L$ the inequality
\[
x\wedge \sim x\leq y\vee \sim y
\]

\noindent holds, then $\mathcal{K}=(L,\vee ,\wedge ,\sim 0,1)$ is called a%
\textit{\ Kleene algebra}. If for $a,b\in L$ there exists a greatest element
$x\in L$ such that $a\wedge x\leq b$, then this $x$ is called \textit{the
relative pseudocomplement of }$a$\textit{\ with respect to} $b$, and it is
denoted by $a\rhd b$. A \textit{quasi-Nelson} algebra is a Kleene algebra $%
\mathcal{K}$ such that $a\rhd (\sim a\vee b)$ exists for all $a$,$b\in L$.
(see e.g. [13]). $a\rhd (\sim a\vee b)$ is denoted simply by $a\rightarrow b$%
. A \textit{Nelson algebra} is an algebra $\mathcal{N}=(A,\vee ,\wedge
,\rightarrow ,\sim 0,1)$ of type (2,2,2,1,0,0), such that $(A,\vee ,\wedge
,\sim 0,1)$ is a quasi-Nelson algebra with $\rightarrow $, and $\rightarrow $
satisfies%
\begin{center}
\hfill $(x\wedge y)\rightarrow z=x\rightarrow (y\rightarrow z)$, for
all $x,y,z \in A$, \hfill (N)
\end{center}

\noindent i.e. the so-called \textit{Nelson-identity}.

Nelson algebras are the algebraic counterparts of the
\textit{constructive logic with strong negation }(see [18, 19]).
Spinks and Veroff proved [22] that to any Nelson algebra
$\mathcal{N}=(A,\vee ,\wedge ,\rightarrow ,\sim 0,1)$
corresponds an integral commutative residuated lattice $\mathcal{L(N)}%
=(A,\vee ,\wedge ,\ast ,\Rightarrow ,0,1)$. For any $ x,y \in A$ the operations $\Rightarrow \;$%
and $\ast $ are defined as follows:%
\[
x\Rightarrow y:=(x\rightarrow y)\wedge (\sim y\rightarrow \;\sim x),
\]%
\[
x\ast y:=\;\sim (x\rightarrow \;\sim y)\vee \sim (y\rightarrow \;\sim x)
\]

\noindent In view of [22] we have $\rceil x:=x\Rightarrow 0=\;\sim
x$, for all $x \in A$, which is an antitone involution. Thus $\rceil
\left( \rceil x\right) =x$, and
applying Theorem 2.40 in [2], we obtain%
\[
x\ast y=\;\rceil (y\Rightarrow \rceil x),
\]

\noindent for all $x,y \in A$, and hence $\Rightarrow $ and the
derived implication of $\mathcal{L(N)}$ coincide. Clearly, the
residuated lattice $\mathcal{L(N)}$ satisfies the
condition (C) and (I3) (see e.g. Corollary 1). Let $%
x^{a}:=x\Rightarrow a$, for all $x,y \in A$. Then for each $a\in L$ the assignment $%
x \mapsto x^{a}$,  $x\in \lbrack a,1]$ is an antitone extensive
mapping, according to Corollary 3. An other important property of
$\mathcal{L(N)}$ is $3$\textit{-potency} (see [22]), which means
that it satisfies the identity:

$x\Rightarrow (x\Rightarrow (x\Rightarrow y))=x\Rightarrow
(x\Rightarrow y)$, for all $x,y \in A$.

\medskip

 Nelson algebras are also fundamental structures in Rough set theory (see [21] or [19]). During the
last decade new approaches have been developed that combine tools of
Fuzzy set theory with that one of Rough set theory, like the
investigations of intuitionistic fuzzy sets, and fuzzy rough sets
(see e.g. [14]). Our expectation is that the algebraic structures
behind these constructions can be reduced to involutive
right-residuated l-groupoids.

\medskip

\subsection{Bounded lattices with sectionally antitone involutions}

\medskip

In this paragraph we are going to show that bounded lattices with
sectionally antitone involutions are common structures equivalent to
involutive right-residuated l-grupoids having the property that
their induced implication $ \Rightarrow $ satisfies a condition
which will be denoted by (I3$^{\ast }$). This will be applied in the
next Section 5.

\medskip

 Let $\mathcal{L}=(L,\vee ,\wedge ,\{^{a}\mid a\in L\},0,1)$ be a
lattice with sectionally antitone mappings $x\mapsto x^{a}$, $x\in
\lbrack a,1]$ and define the operation $x\Rightarrow y:=(x\vee
y)^{y}$, for all $x,y\in L$.
\medskip

\noindent \textbf{Remark 2} Since $(x\vee y)^{y}\geq y$, we have
$\left(
x\Rightarrow y\right) \Rightarrow y=$ $(x\vee y)^{yy}$. Hence the identity $%
\left( x\Rightarrow y\right) \Rightarrow y=x\vee y$, $ x,y \in L $ holds if and only if $%
(x\vee y)^{yy}=x\vee y$, for all $x,y\in L$. Of course, this is
equivalent to the condition that $x^{aa}=x$, for all $a\in L$ and
$x\in \lbrack a,1]$. Therefore, operation $\Rightarrow $ satisfies
the identity
\begin{center}
\hfill $ \left( x\Rightarrow y\right) \Rightarrow y=x\vee y $, for
all $ x,y \in L $
 \hfill (I3*)
\end{center}

\noindent if and only if $\mathcal{L}$ is a lattice with sectionally
antitone involutions. In that case, define $\sim x:=x^{0}$, for all $x\in L$%
. Then $x\mapsto $\/$\sim x$, $x\in L$ is an antitone involution on
the lattice $L$, moreover, $x\Rightarrow 0=x^{0}=$\/$\sim x$, for
all $ x \in L $. \medskip

Since (I3*) implies condition (I3), we can apply Theorem 1 to get:

\bigskip

\noindent \textbf{Theorem 3}
\begin{itemize}
\item[(a)] \emph{Let} $\mathcal{L}=(L,\vee ,\wedge ,\{^{a}\mid a\in L\},0,1)$
\emph{be a bounded lattice with sectionally antitone involutions
}$x\mapsto x^{a},$\emph{\ }$x\in \lbrack a,1]$\emph{. If we define
}$\sim x:=x^{0}$, $x\rightarrow y:=(\thicksim x\vee \thicksim
y)^{\thicksim x}$ \emph{and }$x\odot y:=\thicksim (y\rightarrow
\thicksim x)=\thicksim \lbrack (x\vee \thicksim y)^{\thicksim y}]$,
for all $ x,y \in L $,

\emph{then} $\mathcal{G}(\mathcal{L})=(L,\vee ,\wedge ,\odot
,\rightarrow ,0,1)$ \emph{is an involutive integral right-residuated
l-groupoid with }$\rceil x=\;\thicksim x$\emph{, and its derived implication }%
$\Rightarrow $\emph{\ satisfies }(I3$^{\ast }$)\emph{.}

\item[(b)]  \emph{Let} $\mathcal{G}=(L,\vee ,\wedge ,\odot ,\rightarrow ,0,1)$
\emph{be an involutive integral right-residuated l-groupoid such
that
its derived implication }$\Rightarrow $\emph{\ satisfies condition }(I3$%
^{\ast }$)\emph{, and define }$x^{a}:=x\Rightarrow a$, \emph{for all }$%
a,x\in L$\emph{\ with }$x\geq a$\emph{.}

\emph{Then} $\mathcal{L}(\mathcal{G})=(L,\vee ,\wedge ,\{^{a}\mid
a\in L\},0,1)$ \emph{is bounded lattice with sectionally antitone
involutions }$x\mapsto x^{a},$\emph{\ }$x\in \lbrack a,1]$\emph{, and} $%
x^{0}=x\rightarrow 0$.

\item[(c)] \emph{The correspondence between bounded lattices with sectionally
antitone involutions and involutive integral right-residuated
l-groupoids satisfying }(I3$^{\ast }$)\emph{\ is one-to-one, i.e}. $\mathcal{G%
}(\mathcal{L}(\mathcal{G}))=\mathcal{G}$ \emph{and} $\mathcal{L}(\mathcal{G}(%
\mathcal{L}))=\mathcal{L}$.
\end{itemize}

\medskip

\pf (a) We have to show only that $\mathcal{G}(%
\mathcal{L})=(L,\vee ,\wedge ,\odot ,\rightarrow ,0,1)$ is integral. Since $%
x\Rightarrow 0=x^{0}=$\/$\sim x$ and $x\rightarrow 0=1^{\sim
x}=\;\sim x$ for all $ x \in L $ by
definition, we get $x\Rightarrow 0=x\rightarrow 0$. Hence $\mathcal{G}(%
\mathcal{L})$ is integral, according to Theorem 2(i).

\noindent (b) In view of Theorem 1(b), now it suffices to prove $%
x^{0}=x\rightarrow 0$. Since $\mathcal{G}$ is integral, using the
definition of $\Rightarrow $ and Theorem 2(i) we obtain
$x^{0}=x\Rightarrow 0=x\rightarrow 0$, for all $ x \in L $. (c) is
clear. \ink

\medskip

\noindent \textbf{Proposition 4 }\emph{Let} $\mathcal{G}=(L,\vee
,\wedge ,\odot ,\rightarrow ,0,1)$ \emph{be a right-residuated
l-groupoid. Then the following assertions are equivalent.}

\begin{itemize}
\item[(i)] $(x\Rightarrow y)\Rightarrow y=(y\Rightarrow x)\Rightarrow x$,
\emph{for all} $x,y\in L$\emph{, and} $\mathcal{G}$ \emph{is involutive.}

\item[(ii)] $\Rightarrow $ \emph{satisfies} (I3*)\emph{, and} $\mathcal{G}$
\emph{is involutive.}

\item[(iii)] $\mathcal{G}$ \emph{satisfies the double negation law,
divisibility, and condition }(C).
\end{itemize}

\medskip

\pf (i)$\Rightarrow $(ii). Let $ x,y,z \in L $ arbitrary.
Since\textit{\ }$\mathcal{G} $ is involutive, by Lemma 2 we have
$(x\vee y)\Rightarrow y=x\Rightarrow y$, $1\Rightarrow x=x$, and
$y\leq z$ implies $y\Rightarrow z=1$. Now, using (i) we deduce
(I3*). Indeed,  $(x\Rightarrow y)\Rightarrow y=((x\vee y)\Rightarrow
y)\Rightarrow y=(y\Rightarrow (x\vee y))\Rightarrow (x\vee
y)=1\Rightarrow (x\vee y)=x\vee y$, for all $ x,y \in L $.

\noindent (ii)$\Rightarrow $(iii). Since $\mathcal{G}$ is
involutive, it satisfies the double negation law. Because (I3*)
implies (I3), by Corollary 1 we deduce that $\mathcal{G}$ satisfies
(C) and for any $ x,y \in L $ we have $x\odot y=\rceil (y\rightarrow
\rceil x)$. By using this formula and (I3*) we obtain:

\smallskip
$(x\rightarrow y)\odot x=\rceil (x\rightarrow \rceil
(x\rightarrow y))=\rceil (x\rightarrow \rceil (\rceil y\Rightarrow
\rceil x))=$

$=\rceil ((\rceil y\Rightarrow \rceil x)\Rightarrow \rceil x)=\rceil
(\rceil y\vee \rceil x)=x\wedge y$, for all $ x,y \in L $,

\smallskip
\noindent which proves divisibility.

\noindent (iii)$\Rightarrow $(i). Since $\mathcal{G}$ satisfies (C)
and the double negation law, in view of Proposition 1 it is
involutive, and satisfies $x\odot y=\rceil (y\rightarrow \rceil x)$,
for all $ x,y \in L $. Hence repeating the previous proof we get
$(x\rightarrow y)\odot x=\rceil ((\rceil y\Rightarrow \rceil
x)\Rightarrow \rceil x)$. Now, substituting $x$ by $\rceil x$ and
$y$ by $\rceil y$, for any $ x,y \in L $ we get
\[
\rceil \left( (y\Rightarrow x)\Rightarrow x\right) =(\rceil
x\rightarrow \rceil y)\odot \left( \rceil x\right),
\]

\noindent and then interchanging $x$ and $y$ we obtain:%
\[
\rceil \left( (x\Rightarrow y)\Rightarrow y\right) =(\rceil y\rightarrow
\rceil x)\odot \left( \rceil y\right) .
\]

\noindent Since $(\rceil x\rightarrow \rceil y)\odot \left( \rceil
x\right) =\rceil x\wedge \rceil y=(\rceil y\rightarrow \rceil
x)\odot \left( \rceil y\right) $ by divisibility, we deduce
$(y\Rightarrow x)\Rightarrow x=(x\Rightarrow y)\Rightarrow y$, for
all $ x,y \in L $.  \ink

\medskip

We note that the identity from Proposition 4(i) is called \textit{\L
ukasie\-wicz identity}. Hence we can introduce the following
concept:

\medskip

\noindent \textbf{Definition 2 }If an integral involutive
right-residuated l-groupoid $\mathcal{G}$ satisfies {\L ukasie\-wicz
identity, then we say that $\mathcal{G}$ has \textit{\L ukasiewicz
type. }

\medskip

If $\mathcal{G}$ has \L ukasie\-wicz type, then in view of the proof
of  (ii)$\Rightarrow $(iii) from Proposition 4, $\mathcal{G}$ also
satisfies $x\odot y=\rceil (y\rightarrow \rceil x)$, for all $x,y
\in L$ and (I3).

\smallskip

\section{\L ukasiewicz type right-residuated l-groupoids and basic algebras}

Basic algebras were introduced in [7] and [9] as a common
generalization of MV-algebras and othomodular lattices. The details
of this generalization will be mentioned latter. It is worth
noticing that MV-algebras form an algebraic counterpart of \L
ukasie\-wicz many-valued logic, and othomodular lattices represent
an algebraic framework for certain logical computations motivated by
foundational issues of quantum theory.

\medskip

\noindent \textbf{Definition 3 }By a \textit{basic algebra} is meant
an algebra $\mathcal{A}=(A,\oplus ,\rceil ,0)$ of type $(2,1,0)$
satisfying the following axioms:

\begin{description}
\item[(BA1)] $x\oplus0=x$, for all $ x\in A $

\item[(BA2)] $\rceil \rceil x=x$, for all $ x\in A $

\item[(BA3)] $\rceil (\rceil x\oplus y)\oplus y=\rceil (\rceil y\oplus
x)\oplus x$, for all $ x,y\in A $

\item[(BA4)] $\rceil (\rceil (\rceil (x\oplus y)\oplus y)\oplus z)\oplus
(x\oplus z)=1$, for all $ x,y,z\in A $, where $1:=\rceil 0$.\medskip
\end{description}

\noindent Recall from [7], [8] and [9] that every basic algebra is a
bounded lattice where $x\vee y=\rceil (\rceil x\oplus y)\oplus y$,
$x\wedge y=\rceil (\rceil x\vee \rceil y)$, for all $ x,y\in A $ and
the induced order $\leq $ is given by
\begin{center}
$ x\leq y $ if and only if $ \rceil x\oplus y=1 $.
\end{center}

\noindent Of course, $0\leq x\leq 1,$ for all $x\in A$. In every basic
algebra $\mathcal{A}=(A,\oplus ,\rceil ,0)$ for all $ x,y \in L $ we define the term operations $%
\odot $, $\rightarrow $ and $\Rightarrow $ as follows:%
\begin{center}
$ x\odot y=\rceil (\rceil x\oplus \rceil y)$, $x\rightarrow
y=y\oplus \rceil x$ and
\end{center}
\begin{center}
$ x\Rightarrow y=\rceil x\oplus y $.
\end{center}

\noindent One can observe that $x\Rightarrow 0=\rceil x$, and
$x\Rightarrow y=\rceil y\rightarrow \rceil x$, for all $ x,y \in L
$. The following theorem was established in [9].

\bigskip

\noindent \textbf{Theorem 4}
\begin{itemize}
\item[(i)] \emph{Let} $\mathcal{L}=(L,\vee ,\wedge ,\{^{a}\mid a\in L\},0,1)$
\emph{be a bounded lattice with
sectionally antitone involutions. If we define}%
\begin{center}
$ x\oplus y:=(x^{0}\vee y)^{y}$ \emph{and} $ \rceil x:=x^{0}$,
\emph{for all} $x,y \in L $,
\end{center}

\noindent \emph{then} $\mathcal{A(L)}=(L,\oplus ,\rceil ,0)$ \emph{is a
basic algebra. We have }$x\vee y=\rceil (\rceil x\oplus y)\oplus y$, $%
x\wedge y=\rceil (\rceil x\vee \rceil y)$,
\emph{for all} $x,y \in L $ \emph{and} $x^{a}=\rceil x\oplus a$%
\emph{,} \emph{for }$x\in \lbrack a,1]$.

\item[(ii)] \emph{Let} $\mathcal{A}=(A,\oplus ,\rceil ,0)$ \emph{be a
basic algebra and set }%
\begin{center}
 $x\vee y:=\rceil (\rceil x\oplus y)\oplus y $, $x\wedge y:=\rceil
(\rceil x\vee \rceil y) $, \emph{for all} $x,y \in A $.
\end{center}

\noindent \emph{Define} $x^{a}:=\rceil x\oplus a$, \emph{for all}
$a,x \in A $ \emph{with} $a \leq x$, \emph{and} $1:=\rceil 0$.
\emph{Then} $\mathcal{L(A)}=(A,\vee ,\wedge ,\{^{a}\mid a\in A\},0,1)$ \emph{%
is a bounded lattice with sectionally antitone involutions} $x
\mapsto x^{a}$, $ x \in [a,1]$, \emph{where the
lattice order is given by }$x\leq y$ \emph{iff} $\rceil x\oplus y=1$\emph{%
, and we have} $\rceil x=x^{0}$\emph{,} $x\oplus y:=(x^{0}\vee
y)^{y}$.

\item[(iii)] \emph{The correspondence between bounded lattices with sectionally
antitone involutions and basic algebras thus established is
one-to-one, i.e. }$\mathcal{A}(\mathcal{L(A)})=\mathcal{A}$ \emph{and} $%
\mathcal{L}(\mathcal{A(L)})=\mathcal{L}$\emph{.}
\end{itemize}

\bigskip

\noindent Now, let $\mathcal{A}=(A,\oplus ,\rceil ,0)$ be a basic
algebra and $(A,\vee ,\wedge ,0,1)$ the bounded lattice determined
by $\mathcal{A}$, according to Theorem 4(ii). Then $1:=\rceil 0$,
and in view of Theorem 4(ii)
this is a lattice with sectionally antitone involutions $x\mapsto x^{a}$, $%
x\in \lbrack a,1]$, where $x^{a}:=\rceil x\oplus a$, for all $ a,x \in A $.  In particular, $%
x^{0}=\rceil x$, $x\in A$ determines an involution on the whole
lattice. Further, define
\begin{center}
$x\rightarrow y=(\rceil x\vee \rceil y)^{\rceil x}$ and $x \odot
y=\rceil \lbrack (x\vee \rceil y)^{\rceil y}]$, for all $x,y \in A$.
\end{center}

\noindent Then applying Theorem 3(a) with $\sim x=x^{0}=\rceil x$ we
obtain that $\mathcal{G}(\mathcal{A})=(A,\vee ,\wedge ,\odot
,\rightarrow ,0,1)$ is an involutive integral right-residuated
l-groupoid such that $\Rightarrow $\ satisfies condition (I3$^{\ast
}$). By Proposition 4, the identity
\[
(x\Rightarrow y)\Rightarrow y=(y\Rightarrow x)\Rightarrow x, \mbox{
for all } x,y \in A
\]

\noindent holds, thus $\mathcal{G}(\mathcal{A})$ is of a \L
ukasiewicz type. By Theorem 4(ii), then we obtain $x\oplus
y:=(x^{0}\vee y)^{y}=(\rceil x\vee y)^{y}$. Thus we get $\rceil
(\rceil x\oplus \rceil y)=\rceil \lbrack (x\vee \rceil y)^{\rceil
y}]=x\odot y$, $x\rightarrow y=(\rceil x\vee \rceil y)^{\rceil
x}=y\oplus \rceil x$ and $x\Rightarrow y=\rceil y\rightarrow \rceil
x=\rceil x\oplus y$, for all $ x,y \in A $.

\sk \noindent Conversely, let $\mathcal{G}=(L,\vee ,\wedge ,\odot
,\rightarrow ,0,1)$ be
an involutive right-residuated l-groupoid of \L ukasiewicz type and $%
\Rightarrow $ its derived implication. Then $(x\vee y)\Rightarrow x$
$=y\Rightarrow x$, and $x\odot y=\rceil (y\rightarrow \rceil x)$,
for all $x,y \in L$, in view of
Lemma 2 and Remark 3. For any $a,x\in L$ with $x\geq a$ define $%
x^{a}:=x\Rightarrow a$. Since $\mathcal{G}$ is integral, and by Proposition
4 $\Rightarrow $ satisfies (I3*), we can apply Theorem 3(b) and we get that $%
\mathcal{L}(\mathcal{G})=(L,\vee ,\wedge ,\{^{a}\mid a\in L\},0,1)$ is a
bounded lattice with sectionally antitone involutions\emph{\ }$x\mapsto
x^{a},$\emph{\ }$x\in \lbrack a,1]$, such that $x^{0}=x\rightarrow 0$. Now,
if we define
\begin{center}
$ \rceil x:=x\rightarrow 0 $ and $ x\oplus y:=(x^{0}\vee
y)^{y}=(\rceil x\vee y)^{y} $, for all $x,y \in L $,
\end{center}

\noindent by Theorem 4(i) we obtain a basic algebra $\mathcal{A(G)}%
=(L,\oplus ,\rceil ,0)$, where $x\vee y=\rceil (\rceil x\oplus y)\oplus y$, $%
x\wedge y=\rceil (\rceil x\vee \rceil y)$ and $x^{a}=\rceil x\oplus a$\emph{,%
} for\emph{\ }$x\in \lbrack a,1]$. We get also

\begin{center}
$\rceil \left( \rceil x\odot \rceil y\right) =\rceil (\rceil (\rceil
y\rightarrow \rceil (\rceil x)))=\rceil y\rightarrow \rceil (\rceil
x)=\rceil x\Rightarrow y=$
\end{center}
\begin{center}
$=(\rceil x\vee y)\Rightarrow y=(\rceil x\vee y)^{y}=x\oplus y$, and
\end{center}
\begin{center}
$y\oplus \rceil x=(\rceil y\vee \rceil x)^{\rceil x}=(\rceil y\vee
\rceil x)\Rightarrow \rceil x=\rceil y\Rightarrow \rceil
x=x\rightarrow y$,

\end{center}
\medskip

\noindent for all $x,y \in L$. Now, by using the above computations
we can formulate:

\bigskip

\noindent \textbf{Theorem 5}
\begin{itemize}
\item[(a)] \emph{Let} $\mathcal{A}=(A,\oplus ,\rceil
,0)$ \emph{be a basic algebra. For all} $ x,y \in A $ \emph{define}
\[
x\odot y:=\rceil (\rceil x\oplus \rceil y) \emph{, and }
x\rightarrow y:=y\oplus \rceil x.
\]

\noindent \emph{Set }$x\vee y:=\rceil (\rceil x\oplus y)\oplus y$,
$x\wedge
y:=\rceil (\rceil x\vee \rceil y)$\emph{, and} $1:=\rceil 0$. \emph{Then} $%
\mathcal{G(A)}=(A,\vee ,\wedge ,\odot ,\rightarrow ,0,1)$ \emph{is a
right-residuated l-groupoid of \L ukasiewicz type.}

\item[(b)] \emph{Let} $\mathcal{G}=(A,\vee ,\wedge ,\odot ,\rightarrow ,0,1)$
\emph{be is a right-residuated l-groupoid of \L ukasiewicz type.
Define }$\rceil x:=x\rightarrow 0$ \emph{and} $x\oplus y:=\rceil
\left( \rceil x\odot \rceil y\right) $, \emph{for all} $x,y \in A $.
\emph{Then} $\mathcal{A(G)}=(A,\oplus ,\rceil ,0)$ \emph{is a basic
algebra.}

\item[(c)] \emph{The correspondence between basic algebras and
right-residuated l-groupoids of \L ukasiewicz type thus established is
one-to-one,
i.e. }$\mathcal{A}(\mathcal{G(A)})=\mathcal{A}$ \emph{and} $\mathcal{G}(%
\mathcal{A(G)})=\mathcal{G}$\emph{.}
\end{itemize}

\medskip

\pf Since (a) and (b) follow from the previous computations, we have
to check (c) only. If $\mathcal{A}=(A,\oplus ,\rceil ,0)$ is a basic
algebra, then in $\mathcal{G(A)}$ we have $x\odot y=\rceil (\rceil
x\oplus \rceil y)$, for all $x,y \in A $, and $1=\rceil 0$. Then
$x=1\odot x=\rceil (x\rightarrow \rceil 1)=\rceil (x\rightarrow
\rceil \rceil 0)=\rceil (x\rightarrow 0)$. Thus we get $\rceil
x=\rceil \rceil (x\rightarrow 0)=x\rightarrow 0$, by using (BA2).
This means that $\rceil $ is the same
operation in $\mathcal{A}$ and $\mathcal{A}(\mathcal{G(A)})$. Since in $%
\mathcal{G(A)}$ we have also $\rceil \left( \rceil x\odot \rceil y\right)
=\rceil \rceil (\rceil \rceil x\oplus \rceil \rceil y)=x\oplus y$, in view
of the definition in Theorem 5(b) the operations $\oplus $ in $\mathcal{A}$
and $\mathcal{A}(\mathcal{G(A))}$ coincide. Hence $\mathcal{A}$ and $%
\mathcal{A}(\mathcal{G(A)}$ are the same algebras. The fact that $\mathcal{G}%
(\mathcal{A(G)})=\mathcal{G}$ can be proved similarly. \ink

\bigskip

\noindent The following Corollary is immediate:

\medskip

\noindent \textbf{Corollary 4} \emph{Any} \emph{right-residuated
l-groupoid of \L ukasiewicz type is term equivalent to a basic
algebra. Right-residuated l-groupoids of \L ukasiewicz type form a
variety.}

\medskip

\noindent \textbf{Remark 4} Let $\mathcal{A}=(A,\oplus ,\rceil ,0)$
be a basic algebra, and $x\odot y=\rceil (\rceil x\oplus \rceil y)$,
for all $x,y \in A $. Let us observe that $\odot $ is associative if
and only if $\oplus $ is associative, and $\odot $ is commutative if
and only if $\oplus $ is commutative. Indeed,

\noindent $(x\odot y)\odot z=\rceil \lbrack \rceil (x\odot y)\oplus
\rceil z]=\rceil \lbrack (\rceil x\oplus \rceil y)\oplus \rceil z]$,
and $x\odot \left( y\odot z\right) =\rceil \lbrack \rceil x\oplus
\rceil (y\odot z)]=\rceil \lbrack \rceil x\oplus (\rceil y\oplus
\rceil z)]$. Hence \noindent $(x\odot y)\odot z=x\odot \left( y\odot
z\right)$ if and only if $(\rceil x\oplus \rceil y)\oplus \rceil
z=\rceil x\oplus (\rceil y\oplus \rceil z)$, and this is equivalent
to $(x\oplus y)\oplus z=x\oplus \left( y\oplus z\right) $.

\noindent The proof of the second statement is straightforward.

\bigskip

\textit{Examples}

\medskip

\noindent \textbf{1.} \textit{MV-algebras} form an important
particular case of basic algebras. They can be defined as
associative basic algebras (see e.g. [7]). Since to any basic
algebra corresponds a right-residuated l-groupoid of \L ukasiewicz
type, in view of Remark 4 and Corollary 2, this means that to any
MV-algebra corresponds an integral commutative residuated lattice of
\L ukasiewicz type. We note also that these lattices are always
distributive.

\medskip

\noindent \textbf{2.} \textit{Orthomodular lattices} are usually
defined as bounded orthocomplemented lattices $\mathcal{L=}(L,\vee
,\wedge ,\thicksim ,0,1)$ satisfying the \textit{orthomodular law}
\begin{center}
\hfill $ x\leq y $ implies $ x\vee (\sim x\wedge y)=y $. \hfill
(OML)
\end{center}

\noindent Here $\thicksim $ denotes the \textit{orthocomplementation
operation} on $L$, i.e. $\thicksim $ is an antitone involution such that $%
x\wedge \;\thicksim x=0$, for all $x\in L$.

Define $x^{a}:=\sim x\vee a $, for all $ x,y \in L $. It is known
(see [12] or [4]) that for each $a\in L$ the mapping $x\mapsto x^{a}
$, $x\in \lbrack a,1]$ is an antitone involution on the section
$[a,1]$, moreover $1^{a}=a$. Hence, in view of Theorem 4 (and
Proposition 4), by defining for all $x,y \in L $ the operations
\begin{center}
$x\rightarrow y:=(\thicksim x\vee \thicksim y)^{\thicksim
x}=\;\thicksim (\thicksim x\vee \thicksim y)\vee \thicksim
x=(x\wedge y)\vee \thicksim x$ and
\end{center}
\begin{center}
$x\odot y:=\thicksim \lbrack (x\vee \thicksim y)^{\thicksim
y}]=\;\thicksim \lbrack \thicksim (x\vee \thicksim y)\vee \thicksim
y]=(x\vee \thicksim y)\wedge y$,
\end{center}

\noindent we obtain a right-residuated l-groupoid $\mathcal{G}\left( \mathcal{%
L}\right) =(L,\vee ,\wedge ,\odot ,\rightarrow ,0,1)$ of \L ukasie\-wicz type,
where $\rceil x=\;\thicksim x$. It is easy to check that $\odot $ is not
commutative in general. Therefore, in view of Corollary 2, $\odot $ can not
be even associative.

In [7] was shown that by defining $x\oplus y:=(x\wedge \thicksim
y)\vee y$ for all $x,y\in L$, we obtain a basic algebra
$\mathcal{A}=(L,\oplus ,\rceil ,0)$. It was also proved that basic
algebras arising from orthomodular lattices form a subvariety
characterized by the identity
\begin{center}
 \hfill $ y=y\oplus (x\wedge y) $, for all $x,y\in L$. \hfill (OMI)
\end{center}

\noindent which implies also $x\oplus x=x$, for all $x\in L$. Observe that $%
\mathcal{G}\left( \mathcal{L}\right) $ is just the right-residuated
l-groupoid corresponding to the basic algebra $\mathcal{A}$,
according to
Theorem 5. Now, an easy computation shows that (OMI) is equivalent to $%
\rceil y\rightarrow (x\wedge y)=y$, for all $x,y\in L$. Using the derived implication $%
\Rightarrow $ of $\mathcal{G}\left( \mathcal{L}\right) $, this can
be
reformulated as%
\begin{center}
\hfill $ y=(\rceil x\vee \rceil y)\Rightarrow y $, for all $x,y\in
L$. \hfill (OMI*)
\end{center}

\noindent Hence residuated l-groupoids corresponding to orthomodular
lattices are exactly the right-residuated l-groupoids of \L
ukasiewicz type satisfying (OMI*).

\medskip

\section{Implication reducts of basic algebras}

\medskip

Since the logical connective implication is the most productive one,
because it enables to set up some derivation rules as e.g. Modus
Ponens, we are focused now in a description of implication reducts.

\medskip

Let $\mathcal{A}=(A,\oplus ,\rceil ,0)$ be a basic algebra. For
every $x,y\in A$ define
\[
x\Rightarrow y:=\rceil x\oplus y,
\]

\noindent the so called \textit{implication} in $\mathcal{A}$, and $%
1:=0\Rightarrow 0$. One can easily check that $\Rightarrow $ satisfies the
following identities (see [10]):

\begin{description}
\item[(I0$^{\ast }$)] $x\Rightarrow x=1$, $x\Rightarrow 1=1$, $1\Rightarrow
x=x$, for all $x \in A$;

\item[(I1$^{\ast }$)] $y\Rightarrow (x\Rightarrow y)=1$, for all $x,y\in A$;

\item[(\L )] $(x\Rightarrow y)\Rightarrow y=(y\Rightarrow x)\Rightarrow x$, for all $x,y\in A$;

\item[(I4)] $((x\Rightarrow y)\Rightarrow y)\Rightarrow z)\Rightarrow
(x\Rightarrow z)=1$, for all $x,y\in A$.
\end{description}

\noindent Now, consider the right-residuated l-groupoid $\mathcal{G(A)}%
=(A,\vee ,\wedge ,\odot ,\rightarrow ,0,1)$ which corresponds to the
basic algebra $\mathcal{A}$ by Theorem 5(a). Since $x\rightarrow
y=y\oplus \rceil x $, it is easy to see that $\Rightarrow $
coincides with the so-called
derived implication in $\mathcal{G(A)}$. Since $\mathcal{G(A)}$ is of \L %
ukasiewicz type, in view of Lemma 2 and Proposition 4, for all
$x,y\in A$ the following
assertions also hold true:%
\[
x\leq y\Leftrightarrow x\Rightarrow y=1; \ (x\Rightarrow
y)\Rightarrow y=(x\vee y); \ (x\vee y)\Rightarrow y=x\Rightarrow y.
\]

\noindent Hence the partial order $\leq $ is also determined by $\Rightarrow
$. The fact that $0$ is the least element in $(A,\vee ,\wedge )$, can be
expressed by the law:

\begin{description}
\item[(I5)] $0\Rightarrow x=1$, for all $x\in A$.
\end{description}

\noindent Observe that the previous identities can be inferred from (I0$%
^{\ast }$), (I1$^{\ast }$), (\L ), (I4) and (I5) only, even more, we have
the following

\bigskip

\noindent \textbf{Proposition 5.} \emph{Let} $(A;\Rightarrow ,1)$
\emph{be an algebra of type }$\emph{(2,0)}$\emph{\ satisfying the
identities:}

\begin{itemize}
\item[(i)] $x\Rightarrow x=1$, $x\Rightarrow 1=1$, $1\Rightarrow x=x$, \emph{for all} $x \in A$ \emph{;}

\item[(ii)] $y\Rightarrow (x\Rightarrow y)=1$, \emph{for all} $x,y \in A$\emph{;}

\item[(iii)] $(x\Rightarrow y)\Rightarrow y=(y\Rightarrow x)\Rightarrow x$, \emph{for all} $x,y \in A$;

\item[(iv)] $((x\Rightarrow y)\Rightarrow y)\Rightarrow z)\Rightarrow
(x\Rightarrow z)=1$ \emph{for all} $x,y,z \in A$.
\end{itemize}

\noindent \emph{Define a binary relation }$\leq $\emph{\ on }$A$\emph{\ as
follows}%
\begin{center}
$ x\leq y $ \emph{if and only if} $ x\Rightarrow y=1 $.
\end{center}

\noindent \emph{Then }$\leq $\emph{\ is a partial order on }$A$\emph{, and }$%
(A,\leq )$ \emph{is a join-semilattice with greatest element,}
$1$\emph{\
where }%
\[
x\vee y=(x\Rightarrow y)\Rightarrow y, \mbox{ for all } x,y \in A.
\]

\noindent \emph{Moreover, }$x\leq y$ \emph{implies} $y\Rightarrow z\leq
x\Rightarrow z$\emph{\ and }$\Rightarrow $ \emph{satisfies }%
\[
((x\Rightarrow y)\Rightarrow y)\Rightarrow y=x\Rightarrow y \mbox{
for all } x,y \in A.
\]

\medskip

\pf By (i) the defined relation $\leq $ is reflexive
and $x\leq 1$, for all $x\in A$. Assume $x\leq y$ and $y\leq x$. Then $%
x\Rightarrow y=1$ and $y\Rightarrow x=1$. By (i) and (iii) we
conclude $ y=1\Rightarrow y=(x\Rightarrow y)\Rightarrow
y=(y\Rightarrow x)\Rightarrow x=1\Rightarrow x=x$.

\noindent Let $x\leq y$ and $y\leq z$. Then $x\Rightarrow y=1$ and $%
y\Rightarrow z=1$, and by (iv) we get:

$1=(((x\Rightarrow y)\Rightarrow y)\Rightarrow z)\Rightarrow (x\Rightarrow
z)=((1\Rightarrow y)\Rightarrow z)\Rightarrow (x\Rightarrow z)=$

$=(y\Rightarrow z)\Rightarrow (x\Rightarrow z)=1\Rightarrow (x\Rightarrow
z)=x\Rightarrow z$,

\noindent thus $x\leq z$. Hence $\leq$ is a partial order on $A$ with the
greatest element $1$.

\noindent By (ii) we get $y\leq x\Rightarrow y$, thus also $y\leq
(x\Rightarrow y)\Rightarrow y$ and $x\leq (y\Rightarrow
x)\Rightarrow x $ = $(x\Rightarrow y)\Rightarrow y$, i.e.
$(x\Rightarrow y)\Rightarrow y$ is a common upper bound for $x$ and
$y$.

Next we prove that $a\leq b$ implies $b\Rightarrow c\leq
a\Rightarrow c$. Indeed, $a\leq b$ yields $a\Rightarrow b=1$, and
hence $(b\Rightarrow c)\Rightarrow (a\Rightarrow c)=((1\Rightarrow
b)\Rightarrow c)\Rightarrow (a\Rightarrow c)=(((a\Rightarrow
b)\Rightarrow b)\Rightarrow c)\Rightarrow (a\Rightarrow c)=1$, by
(iv). Hence $b\Rightarrow c\leq a\Rightarrow c$.

Now, if $x,y\leq z$ then $x\Rightarrow y\geq z\Rightarrow y$ and we get also

$(x\Rightarrow y)\Rightarrow y\leq (z\Rightarrow y)\Rightarrow
y=(y\Rightarrow z)\Rightarrow z=1\Rightarrow z=z$,

\noindent proving that $(x\Rightarrow y)\Rightarrow y$ is the least
common upper bound of $x,y$ i.e.

\noindent $(x\Rightarrow y)\Rightarrow y=x\vee y$, for all $ x,y \in
A $. Thus $(A,\leq )$ is a join-semilattice with $1$.

Finally, by using (iii), (ii) and (i), for any $ x,y,z \in A $ we
infer

\noindent $((x\Rightarrow y)\Rightarrow y)\Rightarrow
y=(y\Rightarrow (x\Rightarrow y))\Rightarrow (x\Rightarrow
y)=1\Rightarrow (x\Rightarrow y)=x\Rightarrow y$.
\ink

\bigskip

In what follows, we will consider the algebra $\mathcal{A}%
_{0}=(A,\Rightarrow ,0)$ of type (2,0) which is called an \textit{%
implication reduct} of the basic algebra $\mathcal{A}$. We are going
to show that the basic algebra $(A,\oplus ,\rceil ,0)$ can be
reconstructed from this implication reduct, moreover the following
is true:

\bigskip

\noindent \textbf{Theorem 6.} \emph{Let} $\mathcal{A}_{0}=(A,\Rightarrow ,0)$
\emph{be an algebra of type }$(2,0)$\emph{,} $1:=0\Rightarrow 0$\emph{, such
that} $\Rightarrow $ \emph{satisfies the identities }(i),(ii),(iii),(iv)%
\emph{\ and }(I5)\emph{. Then by defining}%
\begin{center}
\hfill $ \rceil x:=x\Rightarrow 0 $ \emph{and} $ x\oplus y:=\rceil
x\Rightarrow y $, \emph{for all} $ x,y \in A $ \hfill ($ \times $)
\end{center}

\noindent \emph{we obtain a basic algebra} $\mathcal{B}(\mathcal{A}%
_{0})=(A,\oplus ,\rceil ,0)$ \emph{such that the implication in} $\mathcal{B}%
(\mathcal{A}_{0})$ \emph{coincides with} $\Rightarrow $.

\bigskip

\pf In view of Proposition 5, the definition%
\begin{center}
$ x\leq y $ if and only if $ x\Rightarrow y=1 $,
\end{center}

\noindent yields a join-semilattice with greatest element $1$ on the set $A$%
, where $x\vee y=(x\Rightarrow y)\Rightarrow y$, for all $ x,y \in A
$.  In view of (I5), $0$ is the least element of $(A,\leq )$. By
using Proposition 5, we obtain also $\rceil \left( \rceil x\right)
=(x\Rightarrow 0)\Rightarrow 0=x\vee 0=x$, for all $ x \in A $, and
we get that for any $x,y\in A$,
\begin{center}
$ x\leq y $ implies $ \rceil y=y\Rightarrow 0\leq x\Rightarrow
0=\rceil x $.
\end{center}

\noindent This means that the mapping $x\mapsto \rceil x$, $x\in A$
is an antitone involution on $(A,\leq )$, and hence $(A,\leq )$ is a
lattice where $x\wedge y=\rceil (\rceil x\vee \rceil y)$, for all $
x,y \in A $. Since (i),(ii),(iii),(iv)\emph{\ }and\emph{\ }(I5)
together imply the laws (I0),(I1) and (I2) and $(x\Rightarrow
y)\Rightarrow y=x\vee y$, by defining $%
x^{a}:=x\Rightarrow a$ for all $ a,x \in A$, in view of Remark 2, we deduce that the mappings $%
x \mapsto x^{a}$, $x\in \lbrack a,1]$ are antitone involutions on
each section $[a,1]$ of the bounded lattice $(A,\vee ,\wedge )$. In
view of [9] (see Theorem 4), for the operations
\begin{center}
$ x\oplus y:=(x^{0}\vee y)^{y}$ and $ \rceil x:=x^{0} $
\end{center}

\noindent we obtain a basic algebra $(A,\oplus ,\rceil ,0)$ on the
set $A$. Since $x^{0}=x\Rightarrow 0$, $\rceil $ satisfies ($\times
$), and $x\oplus y=(\rceil x\vee y)^{y}=(\rceil x\vee y)\Rightarrow
y=\rceil x\Rightarrow y$,
because (i),(ii),(iii),(iv)\emph{\ }and\emph{\ }(I5)\emph{\ }imply also\emph{%
\ }$(x\vee y)\Rightarrow y=x\Rightarrow y$, for all $ x,y \in A $,
as we pointed out previously.
Finally, the implication in $(A,\oplus ,\rceil ,0)$ is given by the term$%
\rceil x\oplus y$, and $x\oplus y=\rceil x\Rightarrow y$ clearly implies $%
\rceil x\oplus y=x\Rightarrow y$, for all $ x,y \in A $. \ink

\bigskip

We note that Theorem 6 has also a direct proof which does not use
Theorem 4. Observe also, that the conditions (i), (ii), (iii) and
(iv) are in fact the conditions (I0$^{\ast }$), (I1$^{\ast }$), (\L
) and (I4).

\medskip

\section{Congruence properties}

\medskip

When varieties of algebras are studied, we are usually interested in
their congruence properties to reveal their structure.

An algebra $\mathcal{A}=(A,F)$ is said to be \textit{congruence distributive}
whenever its congruence lattice Con$\mathcal{A}$\ is distributive. $\mathcal{%
A}$ is called \textit{congruence permutable}, if $\varphi \circ \theta
=\theta \circ \varphi $ holds for all $\theta ,\varphi \in \;$Con$\mathcal{A}
$. A variety $\mathcal{V}$ of algebras is \textit{arithmetical} if every
algebra $\mathcal{A\in V}$ of it is both congruence distributive and
congruence permutable. An algebra $\mathcal{A}=(A,F)$ is said to be \textit{%
congruence regular} if every congruence $\theta $ of $\mathcal{A}$ is
determined by an arbitrary congruence class $\theta \lbrack a]$ (for $a\in A$%
) of it. Let $c$ be a constant of the algebra $\mathcal{A}$. $\mathcal{A}$
is $c$\textit{-regular} if $\theta \lbrack c]$ = $\varphi \lbrack c]$
implies $\theta =\varphi $, for every $\theta ,\varphi \in \;$Con$A$, and $%
\mathcal{A}$ is called $c$\textit{-locally regular} if for each $\theta
,\varphi \in \;$Con$A$ and any $a\in A$ we have that $\theta \lbrack a]$ = $%
\varphi \lbrack a]$ implies $\theta \lbrack c]$ = $\varphi \lbrack
c]$. It is known that an algebra $\mathcal{A}$ is congruence regular
if and only if it is c-regular and c-locally regular simultaneously
(see [5]). It was proved by B. Cs\'{a}k\'{a}ny [15], that a variety
$\mathcal{V}$ of algebras is congruence c-regular if and only if
there exist binary terms $b_{1},...,b_{n}$ such that $\mathcal{V}$
satisfies the condition
\begin{center}
$ [\ b_{1}(x,y)=c,...,b_{n}(x,y)=c \ ] $ if and only if $x=y$.
\end{center}

\noindent It has been proved in [5] that $\mathcal{V}$ is c-locally
regular if and only if there exist binary terms $p_{1},...,p_{m}$
such that $\mathcal{V}$ satisfies the condition
\begin{center}
$ [\ p_{1}(x,y)=x,...,p_{m}(x,y)=x \ ] $ if and only if $y=c$.
\end{center}

It is known that any right-residuated l-groupoid $\mathcal{G}$ is
congruence 1-regular with the term $b(x,y)=(x\rightarrow y)\wedge
(y\rightarrow x)$ which satisfies $b(x,y)=1 $ if and only if $x=y$.
Clearly, $\mathcal{G}$ is also congruence distributive, because its
reduct to the signature $\{\vee ,\wedge \}$ is a lattice. It is also
known that basic algebras form an arithmetical and congruence
regular variety (see e.g. [7]). Since, in view of Theorem 3,
right-residuated l-groupoids of \L ukasiewicz type are term
equivalent to basic algebras, it follows that they also form an
arithmetical and congruence regular variety. Our last result which
is based on some ideas of [1] shows that some congruence properties
of residuated lattices remain valid in the case of right-residuated
l-groupoids also, although in their case the operation $\odot $ is
neither associative nor integral, in general.

\bigskip

\noindent \textbf{Proposition 6.} \emph{Any} \emph{right-residuated
l-groupoid }$\mathcal{G}=(L,\vee ,\wedge ,\odot ,\rightarrow
,0,1)$\emph{\ is congruence permutable and }$1$\emph{-regular, and
the following hold:}

\begin{itemize}
\item[(a)] \emph{If} $\mathcal{G}$ \emph{satisfies the double negation law,
then it is} $0$\emph{-regular.}

\item[(b)] \emph{If} $\mathcal{G}$ \emph{satisfies divisibility and the
double negation law, then it is congruence regular.}
\end{itemize}

\pf It is well-known that an algebra $\mathcal{A}%
=(A,F) $ is congruence permutable whenever it has a Mal'cev term,
i.e. a term $p(x,y,z)$ satisfying \\
$p(x,y,y)=x$ and $p(x,x,y)=y$, for all $x,y\in A$%
. We can choose the term $p(x,y,z)=\left[ \left( (y\rightarrow z)\wedge
(z\rightarrow y)\right) \odot x\right] \vee \left[ \left( (x\rightarrow
y)\wedge (y\rightarrow x)\right) \odot z\right] $, from [1]. Then $%
p(x,y,y)=x\vee \left[ \left( (x\rightarrow y)\wedge (y\rightarrow
x)\right) \odot y\right] $. Since by Lemma 1(iv) we have $\left(
(x\rightarrow y)\wedge (y\rightarrow x)\right) \odot y\leq
(y\rightarrow x)\odot y\leq x$, we obtain $p(x,y,y)=x$, for all
$x,y\in L$. Similarly we prove $p(x,x,y)=y$, for all $x,y\in L$.

\noindent (a) Let us consider the term $t(x,y)=$ $\left( (x\rightarrow
y)\wedge (y\rightarrow x)\right) \rightarrow 0$. Clearly, $%
t(x,x)=1\rightarrow 0=\rceil 1=0$, according to Lemma 1(vi).
Conversely, if $t(x,y)=0$, then $(x\rightarrow y)\wedge
(y\rightarrow x)=\left( \left( (x\rightarrow y)\wedge (y\rightarrow
x)\right) \rightarrow 0\right)
\rightarrow 0=0\rightarrow 0=1$, by the double negation law. Thus we get $%
x\rightarrow y=1$ and $y\rightarrow x=1$, whence $x\leq y$ and
$y\leq x$, and this implies $x=y$, proving that $\mathcal{G}$ is
$0$-regular.

\noindent (b) Now, in view of (a) and [5], it suffices to prove that
$\mathcal{G}$ is locally $0$-regular. Let $p_{1}(x,y)=$
$(x\rightarrow y)\rightarrow 0$,
and $p_{2}(x,y)=x\vee y$. Then obviously $p_{2}(x,0)=x$, and $%
p_{1}(x,0)=(x\rightarrow 0)\rightarrow 0=x$, for all $ x \in L$.
Conversely, $p_{2}(x,y)=x$
implies $y\leq x$ and $p_{1}(x,y)=x$ yields $(x\rightarrow y)\rightarrow 0=x$%
, whence by double negation we get $x\rightarrow y=x\rightarrow 0$.
Therefore, by using divisibility we obtain: $y=x\wedge
y=(x\rightarrow y)\odot x=(x\rightarrow 0)\odot x=x\wedge 0=0$. This
proves that $\mathcal{G}$ is locally $0$-regular. \ink

\bigskip

\noindent \textbf{Corollary 6.}\emph{\ Let} $\mathcal{V}$ \emph{be a
variety} \emph{consisting of right-residuated l-groupoids}
\emph{satisfying the double negation law and divisibility. Then
}$\mathcal{V}$ \emph{is arithmetical and congruence regular.}


\begin{acknowledgements}
This is a pre-print of an article published as \newline 
I. Chajda, S. Radeleczki, Involutive right-residuated l-groupoids, 
Soft Computing 20 (2015), 119-131, 
doi: 10.1007/s00500-015-1765-7.
The final authenticated version of the article is available online at: 
\newline 
https://link.springer.com/article/10.1007/s00500-015-1765-7.

We also express our thanks to Stephane Foldes for his helpful suggestions.
\end{acknowledgements}



\begin{thebibliography}{10}
\bibitem[1]{1} B\u{e}lohl\'{a}vek R (2003) Some properties of residuated
lattices. Czechoslovak Mathematical Journal 53(128): 161-171


\bibitem[2]{2} B\u{e}lohl\'{a}vek R (2002) Fuzzy relational systems, Foundations
and principles. Kluwer, New-York, ISBN 0-306-46777-1

\bibitem[3]{3} Blyth TS, Janowitz MF (1972) Residuation theory. Int. Ser. of Monographs in Pure and Applied Math.
volume 102. Pergamon Press, Oxford

\bibitem[4]{4} Botur M, Chajda I, Hala\v{s} R (2010) Are basic algebras
residuated structures? Soft Computing 14: 251-255

\bibitem[5]{5} Chajda I (1998) Locally regular varieties. Acta Sci Math
(Szeged) 64: 431-435

\bibitem[6]{6} Chajda I (2003) An extension of relative pseudocomplemention to
non-distributive lattices. Acta Sci Math (Szeged) 69: 491-496

\bibitem[7]{7} Chajda I (2011) Basic algebras and their applications, an
overview. Contributions to General algebra, volume 20.  Verlag J
Heyn. pp 1-10

\bibitem[8]{8} Chajda I, Hala\v{s} R, K\H{u}hr J (2005) Distributive lattices
with sectionally antitone involutions. Acta Sci Math (Szeged) 71:
19-33

\bibitem[9]{9} Chajda I, Hala\v{s} R, K\H{u}hr J (2009) Many-valued quantum
algebras. Algebra Universalis 60: 63-90

\bibitem[10]{10} Chajda I, K\"{u}hr J (2013) Ideals and congruences of basic
algebras. Soft Computing 17: 401-410

\bibitem[11]{11} Chajda I, Radeleczki S (2003) On varieties defined by
pseudocomplemented nondistributive lattices. Publ Math Debrecen
63(4): 737-750

\bibitem[12]{12} Chajda I, Radeleczki S (2014) An approach to orthomodular
lattices via lattices with an antitone involution. Math Slovaca,
accepted for publication

\bibitem[13]{13} Cignoli R (1986) The class of Kleene algebras satisfying an
interpolation property and Nelson algebras. Algebra Universalis 23:
262-292

\bibitem[14]{14} Cornelis C, De Cock M and Radzikowska A M (2007) Vaguely quantified
rough sets. In Rough Sets, Fuzzy Sets, Data Mining and Granular
Computing pp. 87-94.  Springer, Berlin-Heidelberg

\bibitem[15]{15} Cs\'{a}k\'{a}ny B (1970) Characterisation of regular
varieties. Acta Sci Math (Szeged) 31: 187-189

\bibitem[16]{16} Czelakowski J (2008) Additivity of the commutator and
residuation. Reports on Mathematical Logic 43: 109-132.

\bibitem[17]{17} Dilworth R P, Ward M (1939) Residuated lattices. Trans
Amer Math Soc 45: 335-354.

\bibitem[18]{18} J\"{a}rvinen J, Pagliani P, Radeleczki S (2013)
Information completeness in Nelson algebras of rough sets induced by
quasiorders. Studia Logica 101{5}: 1073-1092

\bibitem[19]{19} J�arvinen J and Radeleczki S (2014) Monteiro spaces and rough sets determined by quasiorder relations:
models for Nelson algebras. Fundamenta Informaticae 131 (2):
205-215.

\bibitem[20]{20} Kondo M (2011) Modal operators on commutative residuated
lattices. Math Slovaca 61{1}: 1-14

\bibitem[21]{21} Pagliani P, Chakraborty M (2008) A Geometry of Approximation. Rough
Set Theory: Logic, Algebra and Topology of Conceptual Patterns,
Springer

\bibitem[22]{22} Spinks M, Veroff R (2010) Constructive logic with strong
negation as a substructural logic. J of Logic and Computation 20(4):
761-793

\end{thebibliography}
\end{document}